\documentclass[12pt, authoryear]{elsarticle}              % for format elsarticle
\journal{Advances in Econometrics}            % Name of the Journal, to be used with elsarticle

\usepackage{amssymb,amsfonts,amsmath}
\usepackage{xcolor}
\usepackage{natbib}
\usepackage[top=1.0in, bottom=1.0in, left=1.0in, right=1.0 in]{geometry}   % This is awesome
\usepackage[onehalfspacing]{setspace}
\usepackage[bottom]{footmisc}
\usepackage{indentfirst}
\usepackage{endnotes}
\usepackage{verbatim}                % for block comments
\usepackage{lastpage,placeins}
\usepackage{hyperref}                % for emails. Also, it links all your references...
\usepackage{epsfig,pgf,subfig}       % for including eps,jpeg files. last package is for subfigures
\usepackage{array}
\usepackage{rotating,longtable,booktabs, colortbl}       % needed for rotating the table, makes it better
\usepackage{mdwlist}
\usepackage{enumerate}
\usepackage[shortlabels]{enumitem}
\usepackage{lscape,afterpage}                    % supplies a landscape environment
\usepackage{appendix}
\usepackage{bm}

\usepackage[standard]{ntheorem}
\usepackage[algo2e]{algorithm2e}

\usepackage[T1]{fontenc}
\usepackage{lmodern}

\allowdisplaybreaks[1]                     % page breaks over multi-line equations are allowed, but avoided if possible
\theoremstyle{break}
\theorembodyfont{\normalfont}
\newtheorem{algorithm}[algocf]{Algorithm}

\makeatletter    %changes the category code of '@' character to 11
\def\@author#1{\g@addto@macro\elsauthors{\normalsize%
    \def\baselinestretch{1}%
    \upshape\authorsep#1\unskip\textsuperscript{%
      \ifx\@fnmark\@empty\else\unskip\sep\@fnmark\let\sep=,\fi
      \ifx\@corref\@empty\else\unskip\sep\@corref\let\sep=,\fi
      }%
    \def\authorsep{\unskip,\space}%
    \global\let\@fnmark\@empty
    \global\let\@corref\@empty  %% Added
    \global\let\sep\@empty}%
    \@eadauthor={#1}
}
\makeatother     %reverts this to its original category code of 12.

\begin{document}

\begin{frontmatter} \title{\textbf{Flexible Bayesian Quantile Regression in
Ordinal Models} }

%------------------------------------------------------------------------------
%------------------------------------------------------------------------------
\author{Mohammad Arshad Rahman\corref{cor1}}\ead{marshad@iitk.ac.in}
%
%\fntext[fn1]{Email address: marshad@iitk.ac.in. }
 \address{Department of Economic Sciences, Indian Institute of Technology Kanpur, India.
\\ Office: Room 672, Faculty Building, IIT Kanpur, Kanpur 208016. Phone: +91 512-259-7010. Fax:
+91 512-259-7570.} \cortext[cor1]{Corresponding author}

\author{Shubham Karnawat} \ead{shubkarnawat@gmail.com}
\address{Risk Analyst, Model Risk Management, Credit Suisse, India}

\fntext[fn1]{We thank professor Ivan Jeliazkov (University of California,
Irvine) and the anonymous referee for some valuable comments. Discussion and
suggestions from the participants at the Winter School, Delhi School of
Economics (2016, New Delhi, India), Asian Meeting of the Econometric Society
(2017, Hong Kong), Australasian Meeting of the Econometric Society (2018,
Auckland, New Zealand), and Joint Statistical Meeting (2018, Vancouver,
Canada) are appreciated.}

\fntext[fn2]{Shubham Karnawat worked on this paper while pursuing his M.S.
degree in the Department of Economic Sciences at the Indian Institute of
Technology Kanpur, India.}

\begin{abstract}

This article is motivated by the lack of flexibility in Bayesian quantile
regression for ordinal models where the error follows an asymmetric Laplace
(AL) distribution. The inflexibility arises because the skewness of the
distribution is completely specified when a quantile is chosen. To overcome
this shortcoming, we derive the cumulative distribution function (and the
moment generating function) of the generalized asymmetric Laplace (GAL)
distribution --- a generalization of AL distribution that separates the
skewness from the quantile parameter --- and construct a working likelihood
for the ordinal quantile model. The resulting framework is termed flexible
Bayesian quantile regression for ordinal (FBQROR) models. However, its
estimation is not straightforward. We address estimation issues and propose
an efficient Markov chain Monte Carlo (MCMC) procedure based on Gibbs
sampling and joint Metropolis-Hastings algorithm. The advantages of the
proposed model are demonstrated in multiple simulation studies and
implemented to analyze public opinion on homeownership as the best
long-term investment in the United States (US) following the Great
Recession.

\end{abstract}

%since a single parameter defines both the quantile and the skewness
% Model comparison exhibit the practical utility of the proposed model.

\begin{keyword} Generalized asymmetric Laplace distribution, Gibbs sampling,
Great Recession, homeownership, Markov chain Monte Carlo (MCMC),
Metropolis-Hastings.
\end{keyword}
\end{frontmatter}
%------------------------------------------------------------------------------

%------------------------------------------------------------------------------
\section{Introduction}\label{sec:Intro}
%------------------------------------------------------------------------------

% Layout: 1) Explain when quantile regression is important than a simple
% conditional means model. 2) Discuss why readers should be interested in
% ordinal models, 3) Explain why the inflexibility of the AL is problem and how
% the proposed method improves the previous literature, 4) Finish with your
% roadmap.

Quantile regression, proposed by \citet{Koenker-Basset-1978}, models the
conditional quantiles of the dependent variable as a function of the
covariates. This method is particularly useful if interest lies in the outer
regions of the conditional distribution and/or the data violates the standard
assumptions of mean regression (e.g., presence of heteroscedasticity,
existence of outliers, etc.). Since its introduction, the concept has gained
considerable attention from researchers worldwide and across ideologies.
Within the frequentist econometrics/statistics literature, the advantages of
quantile regression estimators is well studied and the computational
challenges pertaining to optimizing a non-differentiable loss/objective
function have been adequately dealt with. A valuable source is
\citet{KoenkerBook-2005} and references therein. The development of Bayesian
quantile regression faced an impediment since errors in quantile regression
were not assumed to follow any distribution (necessary for writing the
likelihood). About two decades later, \citet{Koenker-Machado-1999} noted that
the quantile loss function appears in the exponent of an asymmetric Laplace
(AL) distribution \citep{Kotz-etal-2001,Yu-Zhang-2005}, thus facilitating the
construction of a parametric likelihood. This distribution was utilized by
\citet{Yu-Moyeed-2001} to propose a Bayesian method for estimating quantile
regression in linear models. The estimation algorithm was further refined in
\citet{Tsionas-2003}, \citet{Reed-Yu-2009} and recently
\citet{Kozumi-Kobayashi-2011} proposed a Gibbs sampling algorithm, where they
exploit the normal-exponential mixture representation of the AL distribution.
The AL likelihood has been utilized to develop algorithms for Bayesian
quantile regression in Tobit models
\citep{Yu-Stander-2007,Kozumi-Kobayashi-2011}, Tobit models with endogenous
covariates \citep{Kobayashi-2017}, censored models \citep{Reich-Smith-2013},
censored dynamic panel data models \citep{Kozumi-Kobayashi-2012}, count data
models \citep{Lee-Neocleous-2010} and mixed-effect or longitudinal data
models \citep{Geraci-Bottai-2007,Luo-Lian-Tian-2012}.

Quantile regression in ordinal models is different since the dependent
variable takes discrete and ordered values (which has no cardinal
interpretation), and does not yield continuous quantiles. Ordinal outcomes
typically arise as response to surveys, and applications are common in
economics, finance, marketing, and the social sciences. Similar to the
continuous case, interest in ordinal quantile regression is aimed to provide
a much richer view of the heterogeneous effect of the covariates on the
outcomes. However, estimation is more challenging. A frequentist approach
using simulated annealing was proposed in \citet{Zhou-2010}. Bayesian
estimation of ordinal quantile regression was introduced in
\citet{Rahman-2016} and extended to longitudinal data models in
\citet{Alhamzawi-Ali-2018}. A special case of ordinal model is the binary
model, where the outcome variable is dichotomous (i.e. takes only two values,
typically coded as 1 for `success' and 0 for `failure'). Bayesian quantile
regression in binary models was proposed in \citet{Benoit-Poel-2010} and
employed to study the mode of transportation to work.
\citet{Rahman-Vossmeyer-2019} extended Bayesian quantile regression to binary
longitudinal outcomes and proposed an efficient Markov chain Monte Carlo
(MCMC) algorithm for its estimation. The Bayesian ordinal quantile regression
model has been utilized in a wide variety of studies including evaluation of
credit risk \citep{Migueis-etal-2013}, educational attainment
\citep{Rahman-2016}, public opinion on tax policy \citep{Rahman-2016}, public
opinion on nuclear power plants operation \citep{Omata-etal-2017}, and
illness severity \citep{Alhamzawi-Ali-2018}.

The list of papers on Bayesian quantile regression mentioned above, although
incomplete, clearly affirm that the AL distribution has played a crucial role
in the development of Bayesian quantile regression. However, the AL
distribution poses a critical limitation since a single parameter defines
both the quantile and the skewness of the distribution. In addition, the mode
of the distribution is always fixed at the location parameter value for all
quantiles. To overcome these drawbacks, \citet{Yan-Kottas-2017} proposed the
probability density function (\emph{pdf}) of the generalized asymmetric
Laplace (GAL) distribution by introducing a shape parameter into the mean of
the normal kernel in the AL mixture representation. The GAL distribution uses
different parameters for quantile and skewness, and thus adds much needed
flexibility for Bayesian quantile regression. They utilized the GAL
\emph{pdf} and proposed algorithms for Bayesian quantile estimation of linear
models, Tobit models and regularized quantile regression.

In this paper, we present a derivation of the GAL \emph{pdf} from the mixture
representation and both introduce and derive the cumulative distribution
function (\emph{cdf}) and the moment generating function (\emph{mgf}) of the
GAL distribution. The GAL density and the GAL \emph{cdf} are utilized to
introduce an estimation method for the flexible Bayesian quantile regression
in ordinal (FBQROR) models. Estimation of ordinal models, unlike linear
models, is more challenging since there are identification restrictions and
sampling of cut-points have to satisfy the ordering constraints. Moreover,
through careful transformation of the mixture variables and joint sampling of
the scale and shape parameters, we are able to achieve low autocorrelation in
our Markov chain Monte Carlo (MCMC) draws. This result is a substantial
improvement compared to the extremely high autocorrelation reported in
\citet{Yan-Kottas-2017}. Our sampling scheme can therefore improve the
algorithm for Bayesian quantile regression in linear, Tobit and regularized
regression models as presented in \citet{Yan-Kottas-2017}.

We illustrate the proposed methodology in two simulation studies where the
errors are generated from a symmetric (logistic) distribution and an
asymmetric (chi-square) distribution. The results show that the FBQROR model
can maintain the actual skewness of the data across all considered quantiles.
Furthermore, the FBQROR models can provide better model fit compared to the
fit obtained from Bayesian quantile regression in ordinal (BQROR) models
assuming an AL distribution. Finally, we implement our FBQROR model in an
application related to the recent housing crisis and the Great Recession (Dec
2007 - Jun 2009). Specifically, we analyze how various socioeconomic \&
demographic factors and exposure to financial distress are associated with
differences in views on the financial benefits of homeownership following the
Great Recession. The results increase our understanding and offer new
insights which may be important for policymakers and researchers interested
in the US housing market.

The remainder of the paper is organized as follows. Section~\ref{sec:GAL}
presents some fundamental properties of the GAL distribution.
Section~\ref{sec:fbqror} presents the FBQROR model and its estimation
procedure. Section~\ref{sec:simStudies} illustrates the algorithm in two
simulation studies and Section~\ref{sec:application} implements the algorithm
to examine US public opinion on homeownership. Section~\ref{sec:conclusion}
presents some concluding remarks.

%------------------------------------------------------------------------------
\section{The GAL Distribution}\label{sec:GAL}
%------------------------------------------------------------------------------

The GAL distribution is obtained by introducing a shape parameter into the
mean of the normal kernel in the normal-exponential mixture representation of
the AL distribution and mixing with respect to a half-normal distribution.
This hierarchical representation allows the skewness and mode to vary for a
given quantile/percentile and hence provides the much needed flexibility for
Bayesian quantile regression.

Suppose $Y$ is a random variable that has the following mixture
representation,
%-------------------------------
\begin{equation}
Y = \mu + \sigma A W + \sigma \alpha S +  \sigma [B W ]^{\frac{1}{2}} U,
\label{eq:GALmixture}
\end{equation}
%-------------------------------
where $W \sim \mathcal{E}(1)$, $S\sim N^{+}(0,1)$, $U\sim N(0,1)$, $A\equiv
A(p)=\frac{1-2p}{p(1-p)}$ and $B\equiv B(p)=\frac{2}{p(1-p)}$. Here,
$\mathcal{E}, N^{+}$ and $N$ denote exponential, half-normal and normal
distributions, respectively. Then, $Y$ follows a GAL distribution denoted $Y
\sim GAL(\mu, \sigma, p, \alpha)$ and has the \emph{pdf},
%-------------------------------
\begin{equation}
\begin{split}
f(y|\theta)  & = \frac{2p(1-p)}{\sigma} \Bigg( \bigg[
\Phi\left(\frac{y^{\ast}}{\alpha} - \alpha p_{\alpha_{-}}\right)
- \Phi\left(- \alpha p_{\alpha_{-}} \right)    \bigg]
\exp\bigg\{ - y^{\ast} p_{\alpha_{-}}
+ \frac{1}{2} \big(\alpha p_{\alpha_{-}}\big)^{2}
\bigg\} \\
& \times I\left( \frac{y^{\ast}}{\alpha} > 0\right)  +
\Phi\left(\alpha p_{\alpha_{+}} -
\frac{y^{\ast}}{\alpha} I\left( \frac{y^{\ast}}{\alpha} > 0\right) \right)
\exp\bigg\{ - y^{\ast} p_{\alpha_{+}} +
\frac{1}{2} \big(\alpha p_{\alpha_{+}} \big)^{2}
\bigg\}
\Bigg), \label{eq:GALpdf}
\end{split}
\end{equation}
%-------------------------------
where $\theta=(\mu,\sigma,p,\alpha)$, $y^{\ast} = (y-\mu)/\sigma$, $\mu$ is
the location parameter, $\sigma$ is the scale parameter, $\alpha$ is the
shape parameter, $p_{\alpha_{+}} = p -I(\alpha>0)$ and $p_{\alpha_{-}} = p
-I(\alpha<0)$ with $p \in (0,1)$. The derivation of the GAL \emph{pdf} from
the hierarchical representation is presented in \ref{app:GALpdf} and largely
follows the notations used in \citet{Yan-Kottas-2017}. Note that when
$\alpha=0$, the GAL \emph{pdf} reduces to the \emph{pdf} of an AL
distribution.

We explore the GAL distribution in greater detail and propose the \emph{cdf}
and \emph{mgf} of the GAL distribution. The \emph{cdf} denoted by $F$ can be
compactly written as,
%---------------------------------
\begin{equation}
\begin{split}
F(y|\theta) = \bigg(1 - 2 \Phi\left(-\frac{y^{\ast}}{|\alpha|} \right)
- \frac{2p(1-p)}{p_{\alpha_{-}}} \exp\left\{- y^{\ast} p_{\alpha_{-}}
+ \frac{1}{2} \alpha^{2} p_{\alpha_{-}}^{2}  \right\}
\Big[\Phi\left(\frac{y^{\ast}}{\alpha} - \alpha p_{\alpha_{-}} \right) \\
%----
- \, \Phi(-\alpha p_{\alpha_{-}}) \Big] \bigg) I\left( \frac{y^{\ast}}{\alpha}>0
\right) + I(\alpha <0) - \frac{2p(1-p)}{p_{\alpha_{+}}}
\exp\left\{-y^{\ast} p_{\alpha_{+}} + \frac{1}{2} \alpha^{2} p_{\alpha_{+}}^{2}
\right\} \\
%-----
\times \, \Phi\left[\alpha p_{\alpha_{+}} - \frac{y^{\ast}}{\alpha}
I\left( \frac{y^{\ast}}{\alpha}>0\right)   \right],
\end{split} \label{eq:GALcdf}
\end{equation}
%---------------------------------
and the \emph{mgf} denoted by $M_{Y}(t)$ has the following expression,
%---------------------------------
\begin{equation}
M_{Y}(t) = 2p(1-p) \bigg[\frac{ (p_{\alpha_{+}} - p_{\alpha_{-}}) }
{(p_{\alpha_{-}} -\sigma t)(p_{\alpha_{+}} -\sigma t)} \bigg]
\exp\Big\{ \mu t + \frac{1}{2} \alpha^2 \sigma^2 t^2 \Big\} \;
\Phi\big(|\alpha| \sigma t \big). \label{eq:GALmgf}
\end{equation}
%---------------------------------
Both the \emph{cdf} and \emph{mgf} have been derived and presented in
\ref{app:GALcdf} and \ref{app:GALmgf}, respectively. In addition,
\ref{app:GALmgf} utilizes the \emph{mgf} \eqref{eq:GALmgf} to derive the
mean, variance and skewness of the distribution. These distributional
characteristics are extremely important for better understanding of the GAL
distribution and for further development of flexible Bayesian quantile
regression.

However, the GAL density given by equation~\eqref{eq:GALpdf} has the
limitation that the parameter $p$ no longer corresponds to the cumulative
probability at the quantile for $\alpha \neq 0$. We let $\gamma =
[I(\alpha>0) - p]|\alpha|$ and re-express the mixture representation
\eqref{eq:GALmixture} as follows,
%-------------------------------
\begin{equation}
Y = \mu + \sigma A W + \sigma C |\gamma| S +  \sigma [BW]^{\frac{1}{2}} U,
\label{eq:GALmixture2}
\end{equation}
%-------------------------------
where $C=[I(\gamma>0)-p]^{-1}$. This re-parametrization yields the
quantile-fixed GAL distribution that has the following \emph{pdf}:
%-------------------------------
\begin{align}
%\begin{split}
f_{p_{0}}(y|\eta)  & = \frac{2p(1-p)}{\sigma} \Bigg( \bigg[
\Phi\left(-y^{\ast}\, \frac{p_{\gamma_{+}}}{|\gamma|} + \frac{p_{\gamma_{-}}}
{p_{\gamma_{+}}} |\gamma|\right) -
\Phi\left(\frac{p_{\gamma_{-}}}{p_{\gamma_{+}}}|\gamma| \right) \bigg]
\exp\bigg\{ - y^{\ast} p_{\gamma_{-}} + \frac{\gamma^2}{2}
\bigg( \frac{p_{\gamma_{-}}}{p_{\gamma_{+}}} \bigg)^{2}
\bigg\} \nonumber \\
& \times I\left( \frac{y^{\ast}}{\gamma} > 0\right)  +
\Phi\left(- |\gamma| + y^{\ast}\, \frac{ p_{\gamma_{+}}}{|\gamma|}
I\left( \frac{y^{\ast}}{\gamma} > 0\right) \right)
\exp\bigg\{ - y^{\ast} p_{\gamma_{+}} +
\frac{\gamma^2}{2} \bigg\}
\Bigg), \label{eq:GALpdf2}
%\end{split}
\end{align}
%-------------------------------
where $\eta=(\mu,\sigma,\gamma)$, $p\equiv p(\gamma,p_{0})=I(\gamma<0) +
[p_{0} - I(\gamma<0)]/g(\gamma)$, $p_{\gamma_{+}} = p - I(\gamma > 0)$ and
$p_{\gamma_{-}} = p - I(\gamma < 0)$. The function $g(\gamma) = 2
\Phi(-|\gamma|) \exp(\gamma^{2}/2)$ and $\gamma \in (L,U)$, where $L$ is the
negative square root of $g(\gamma)=1-p_{0}$ and $U$ is the positive square
root of $g(\gamma)=p_{0}$. The term ``quantile-fixed'' suggests that
integration of GAL \emph{pdf} \eqref{eq:GALpdf2} to the upper limit $\mu$
equals $p_{0}$, so for regression purpose we can fix the quantile. The
\emph{cdf} for the quantile-fixed GAL density \eqref{eq:GALpdf2} can be
analogously derived as in \ref{app:GALcdf} to yield the following expression,
%-------------------------------
%---------------------------------
\begin{equation}
\begin{split}
F_{p_{0}}(y|\eta) = \Bigg( 1 - 2\Phi\left(y^{\ast} \,\frac{p_{\gamma_{+}}}
{\gamma} \right)  + 2 p_{\gamma_{+}} \exp\bigg\{-y^{\ast} p_{\gamma_{-}}
+ \frac{\gamma^2}{2} \bigg( \frac{p_{\gamma_{-}}}{p_{\gamma_{+}}} \bigg)^{2}
\bigg\}  \bigg[
\Phi\left(- y^{\ast} \,\frac{ p_{\gamma_{+}}}{|\gamma|} + \frac{p_{\gamma_{-}}}
{p_{\gamma_{+}}} |\gamma|\right)  \notag \\
%-----
- \Phi\left(\frac{p_{\gamma_{-}}}{p_{\gamma_{+}}}|\gamma| \right) \bigg]
\Bigg) I\left( \frac{y^{\ast}}{\gamma} >0\right) + I(\gamma<0)
+ 2 p_{\gamma_{-}} \exp\bigg\{-y^{\ast} p_{\gamma_{+}} +
\frac{\gamma^{2}}{2} \bigg\}  \notag \\
%-----
\times \Phi\left(-|\gamma| + y^{\ast} \, \frac{p_{\gamma_{+}}}
{|\gamma|} I\left(\frac{y^{\ast}}{\gamma}>0 \right) \right).
\end{split}
%---------------------------------
\label{eq:GALcdf2}
\end{equation}
%-------------------------------
The quantile-fixed \emph{cdf} \eqref{eq:GALcdf2} is required for constructing
the likelihood of the FBQROR model and plays a critical role in the MCMC
sampling of the scale parameter, shape parameter and cut-points or
thresholds.

%----------------------------  Figure 1 ----------------------------------------
\begin{figure*}[!t]
  \centerline{
    \mbox{\includegraphics[width=7.75in, height=2.8in]{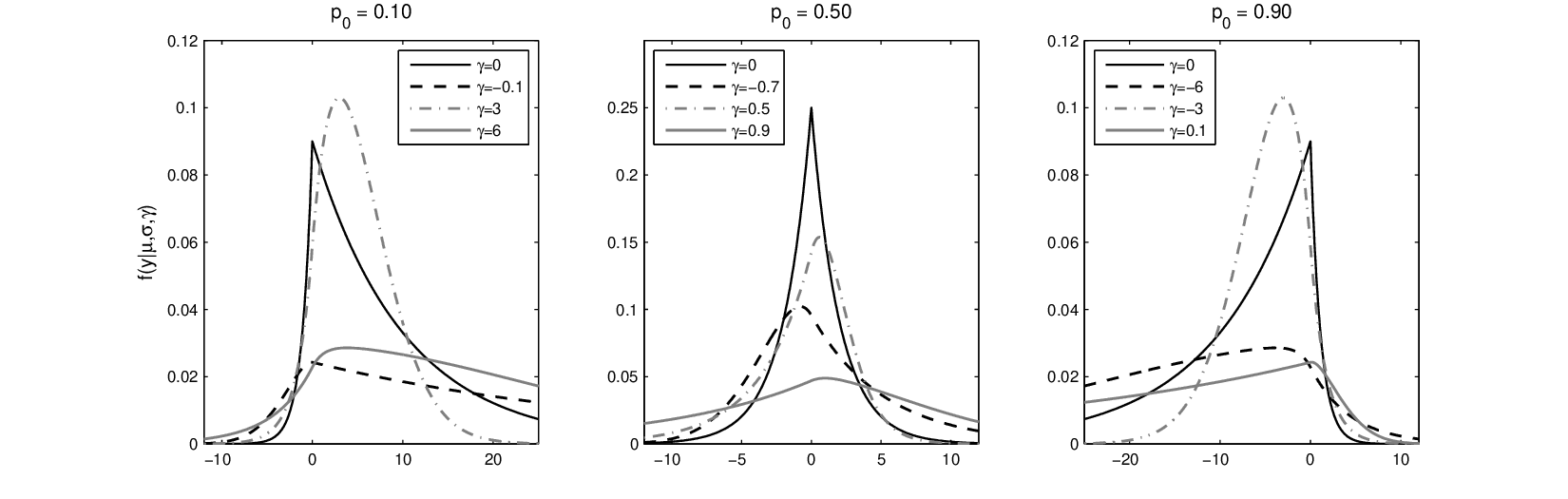}}
    }
\caption{Probability density plots of the AL ($\gamma=0$) and the GAL
($\gamma \neq 0$) distributions.} \label{fig:galpdf}
\end{figure*}
%----------------------------------------------------------------------------------------

To better discern the GAL distribution, Figure~\ref{fig:galpdf} presents a
graphical comparison of the quantile-fixed GAL and AL \emph{pdf}'s for three
different quantiles. We observe that the GAL distribution, unlike the AL
distribution, allows the mode to vary rather than being fixed at $\mu=0$ at
all quantiles. Besides, the GAL distribution can be positive or negatively
skewed at all quantiles depending on the value of $\gamma$. For example, at
the median $p_{0}=0.50$ the GAL distribution is positively skewed for
$\gamma<0$ and negatively skewed for $\gamma>0$. Also, the GAL distribution
can have tails that are heavier or narrower than the AL distribution. These
characteristics make the GAL distribution more flexible than the AL
distribution.

%------------------------------------------------------------------------------
\section{The FBQROR Model}\label{sec:fbqror}
%------------------------------------------------------------------------------

Ordinal models arise when the dependent (response) variable is discrete and
outcomes are inherently ordered or ranked such that the scores assigned to
outcomes have an ordinal meaning, but no cardinal interpretation
\citep{Johnson-Albert-2000,Jeliazkov-Rahman-2012}. For example, in a survey
on public opinion to allow more offshore drilling, responses may be recorded
as follows: 1 for `strongly oppose', 2 for `somewhat oppose', 3 for `somewhat
support' and 4 for `strongly support' \citep{Mukherjee-Rahman-2016}. These
responses have ordinal meaning but no cardinal interpretation. Therefore, one
cannot say that a score of 4 implies four times more support compared to a
score of 1.

We adopt the latent variable approach and represent the FBQROR model using a
continuous latent random variable $z_{i}$ expressed as a function of
covariates and error as,
%------------------------------
\begin{equation}
z_{i} = x'_{i} \beta  + \epsilon_{i}, \hspace{0.75in} \forall \; i=1, \cdots, n,
\label{eq:model1}
\end{equation}
%------------------------------
where $x_{i}$ is a $k \times 1$ vector of covariates, $\beta$ is a $k \times
1$ vector of unknown parameters at the $p_{0}$-th quantile, $\epsilon_{i}$
follows a GAL distribution, i.e., $\epsilon_{i} \sim GAL(0,\sigma,\gamma)$
and $n$ denotes the number of observations. Note that we have suppressed the
dependence of parameters on $p_{0}$ for notational simplicity. The variable
$z_{i}$ is unobserved and relates to the observed discrete response $y_{i}$,
which has \textit{J} categories or outcomes, via the cut-point vector $\xi$
as follows:
%-------------------------------------
\begin{equation}
\xi_{j-1} < z_{i} \le \xi_{j} \; \Rightarrow \;
\emph{$y_{i}$ = j}, \hspace{0.75in}
\forall \; i=1,\cdots, n; \; j=1,\cdots, J,
\label{eq:cutpoints}
\end{equation}
%-------------------------------------
where $\xi_{0}=-\infty$ and $\xi_{J}=\infty$. In addition, $\xi_{1}$ is
typically set to 0, which anchors the location of the distribution required
for parameter identification (see \citealp{Jeliazkov-etal-2008}). Given the
data vector $y$ = $(y_{1}, \cdots, y_{n})'$, the likelihood for the model
expressed as a function of unknown parameters $(\beta, \sigma,\gamma,\xi)$
can be written as,
%-------------------------------------
\begin{equation}
\begin{split}
f(\beta,\sigma,\gamma, \xi; y)
    & = \prod_{i=1}^{n} \prod_{j=1}^{J} P(y_{i} = j | \beta,\sigma,\xi,
        \gamma)^{ I(y_{i} = j)}  \\
    & = \prod_{i=1}^{n}  \prod_{j=1}^{J}
        \bigg[F_{p_{0}}\left(\frac{\xi_{j} - x'_{i}\beta}{\sigma}\right)
        - F_{p_{0}}\left(\frac{\xi_{j-1} - x'_{i}\beta}{\sigma}\right)
        \bigg]^{ I(y_{i} = j)}
            \label{eq:fulllikelihood}
\end{split}
\end{equation}
%-------------------------------------
where, $F_{p_{0}}(\cdot) \equiv F(\cdot|0,1,\gamma)$ denotes the \emph{cdf}
of the GAL distribution and $I(y_{i}=j)$ is an indicator function, which
equals 1 if $y_{i}=j$ and 0 otherwise.

Working directly with the GAL distribution is difficult, so we replace the
error term with its mixture representation \eqref{eq:GALmixture2} and rewrite
the FBQROR model as follows:
%---------------------------
\begin{equation}
z_{i} = x'_{i} \beta   + \sigma A w_{i} + \sigma C |\gamma| s_{i}
        + \sigma  \sqrt{B w_{i}} \,u_{i},
        \hspace{0.5in} \forall \; i=1, \cdots, n.
\label{eq:model2}
\end{equation}
%---------------------------
The above formulation \eqref{eq:model2} implies that the latent variable
$z_{i}|\beta,w_{i},s_{i},\sigma,\gamma \sim N( x'_{i}\beta + \sigma C
|\gamma| s_{i} + \sigma A w_{i}, \sigma^{2} B w_{i})$. However, the presence
of the scale parameter $\sigma$ in the conditional mean is not conducive to
the construction of MCMC algorithm \citep{Kozumi-Kobayashi-2011}. We
reparameterize and write the model as,
%---------------------------
\begin{equation}
z_{i} = x'_{i} \beta  + A \nu_{i} + C |\gamma| h_{i}
        + \sqrt{\sigma B \nu_{i}} \,u_{i},
        \hspace{0.5in} \forall \; i=1, \cdots, n,
\label{eq:model3}
\end{equation}
%---------------------------
where $h_{i}=\sigma s_{i}$ and $\nu_{i} = \sigma w_{i}$, which in turn
implies that $h_{i}\sim N^{+}(0,\sigma^{2})$ and $\nu_{i} \sim
\mathcal{E}(\sigma)$ for $i=1,\cdots,n$. Both reformulations are necessary
for computational efficiency and low autocorrelation in MCMC draws. Note that
the first reparameterization was not utilized in \citet{Yan-Kottas-2017} and
hence our approach provides a better alternative to estimating linear, Tobit
and regularized lasso quantile regression.

Ordinal models present two additional challenges: location and scale
restrictions for identification of the parameters and ordering constraints in
sampling of cut-points $\xi$ \citep[see][]{Jeliazkov-etal-2008,Rahman-2016}.
In the FBQROR model, both location and scale restrictions are enforced by
fixing two cut-points since the variance of a GAL distribution is not fixed
due to its dependence on $\alpha$ even if we set $\sigma=1$, as shown in
Theorem~4 in \ref{app:GAL}. The ordering constraint is resolved by using the
following logarithmic transformation,
%---------------------------------------
\begin{equation}
\delta_{j} = \ln ( \xi_{j+2} - \xi_{j+1} ), \qquad 1 \le j \le J-3.
\label{eq:logtransformation}
\end{equation}
%-----------------------------------------
The original cut-points are then obtained using
equation~\eqref{eq:logtransformation} by one-to-one mapping between $\delta =
(\delta_{1}, \cdots, \delta_{J-3})'$ and $\xi=(\xi_{3}, \cdots, \xi_{J-1})'$,
where $\xi_{2}$ is fixed at some constant $c$, and recall that $\xi_{0} =
-\infty$, $\xi_{1} = 0$ and $\xi_{J} = \infty$.

%----------------------------  Algorithm 1 -------------------------------------
\begin{table*}
\begin{algorithm}[Sampling in FBQROR Model] \label{alg:algorithm1}
\rule{\textwidth}{0.5pt} \small{
\begin{enumerate}[(1)]
%-------------------------------------------------
\item Sample $\beta|z, \nu, h, \sigma, \gamma \sim
    N(\tilde{\beta},\tilde{B})$, where
\begin{equation*}
\tilde{B}^{-1} = \Bigg( B_{0}^{-1} + \sum_{i=1}^{n}
\frac{x_{i}x_{i}^{\prime}}{\sigma B\nu_{i}} \Bigg) \hspace{0.15in}
\text{and} \hspace{0.15in}
\tilde{\beta} = \tilde{B} \Bigg(\sum_{i=1}^{n}
\frac{x_{i}(z_{i}-A\nu_{i}- C |\gamma |h_{i})}{\sigma B\nu_{i}}
+ B_{0}^{-1}\beta_{0} \Bigg).
\end{equation*}
%-------------------------------------------------
\item Sample $(\sigma,\gamma)$ marginally of $(z,\nu,h)$ using a joint
    random-walk MH algorithm. The proposed values $(\sigma',\gamma')$ are
    generated from a truncated bivariate normal distribution
    $TBN_{(0,\infty) \times (L,U)}\big( (\sigma_{c},\gamma_{c}), \iota_{1}^{2}
    \hat{D}_{1}\big)$, where $(\sigma_{c},\gamma_{c})$ denote the current values,
    $\iota_{1}$ denotes the tuning factor and $\hat{D}_{1}$ is the negative inverse
    of the Hessian obtained by maximizing the log-likelihood~\eqref{eq:fulllikelihood}
    with respect to $(\sigma, \gamma)$. The proposed draws are accepted with
    MH probability,
    %--------------------------
    \begin{equation*}
    \alpha_{MH}(\sigma_{c},\gamma_{c};\sigma',\gamma') = \min
    \bigg\{0,\ln\bigg[\frac{f(y| \beta,\sigma',\gamma',\delta) \,
    \pi(\beta,\sigma',\gamma',\delta)}{f(y|
    \beta,\sigma_{c},\gamma_{c},\delta) \, \pi(\beta,\sigma_{c},\gamma_{c},\delta)} \;
    \frac{\pi(\sigma_{c},\gamma_{c}| (\sigma',\gamma'), \iota_{1}^{2} \hat{D}_{1})}
    {\pi(\sigma',\gamma'| (\sigma_{c},\gamma_{c}), \iota_{1}^{2} \hat{D}_{1})}
    \bigg] \bigg\},
    \end{equation*}
    %--------------------------
    else, repeat $(\sigma_{c},\gamma_{c})$ in the next MCMC iteration.
    Here, $f(\cdot)$ represents the full likelihood
    \eqref{eq:fulllikelihood} obtained as the difference of \emph{cdf},
    $\pi(\beta,\sigma,\delta,\gamma)$ denotes the
    prior distributions~\eqref{eq:priors}, and
    $\pi(\sigma_{c},\gamma_{c}| (\sigma',\gamma'), \iota_{1}^{2} \hat{D}_{1})$ stands
    for the truncated bivariate normal probability with mean $(\sigma',\gamma')$
    and covariance $\iota_{1}^{2} \hat{D}_{1}$. The term
    $\pi(\sigma',\gamma'| (\sigma_{c},\gamma_{c}), \iota_{1}^{2} \hat{D}_{1})$ has an
    analogous interpretation.
%-------------------------------------------------
\item Sample $\nu_{i}|z_{i}, \beta, h, \sigma, \gamma \sim GIG(0.5,
    a_{i},b)$, for $i=1,\ldots,n$, where
    \begin{equation*}
    a_{i} = \frac{(z_{i}-x_{i}^{\prime}\beta- C |\gamma |
             h_{i})^2}{\sigma B} \hspace{0.15in} \text{and} \hspace{0.15in}
    b = \bigg(\frac{A^2}{\sigma B} + \frac{2}{\sigma} \bigg).
    \end{equation*}
%-------------------------------------------------
%-------------------------------------------------
\item Sample $h_{i}|z_{i},\beta,\nu_{i},\sigma,\gamma \sim N^{+}
        ( \mu_{h_{i}}, \sigma_{h_{i}}^{2} )$ for $i=1,\ldots,n$, where
    %--------------------------
    \begin{equation*}
    (\sigma_{h_{i}}^{2})^{-1} = \bigg(\frac{1}{\sigma^2} + \frac{C^{2}
    \gamma^{2}}{\sigma  B \nu_{i}}  \bigg)
             \hspace{0.25in} \mathrm{and} \hspace{0.25in}
    \mu_{h_{i}} = \sigma_{h_{i}}^{2} \bigg( \frac{ C|\gamma|(z_{i}- x'_{i}\beta
    - A\nu_{i})}{\sigma B \nu_{i}}  \bigg).
    \end{equation*}
    %--------------------------
%-------------------------------------------------
\item Sample $\delta|\beta,\sigma,\gamma,y$ marginally of $(z,\nu,h)$ using a
    random-walk MH step. The proposed value $\delta'$ is generated as
    $\delta' = \delta_{c} + u $, where $u \sim N(0_{J-3}, \iota_{2}^2 \hat{D}_{2})$,
    $\iota_{2}$ is a tuning parameter and $\hat{D}_{2}$ is analogous to $\hat{D}_{1}$.
    Accept $\delta'$ with MH probability,
%--------------------------
\begin{equation*}
\alpha_{MH}(\delta_{c},\delta') = \min \bigg\{0,\ln\bigg[\frac{f(y|
\beta,\sigma,\gamma,\delta') \, \pi(\beta,\sigma,\gamma,\delta')}{f(y|
\beta,\sigma,\gamma,\delta_{c}) \, \pi(\beta,\sigma,\gamma,\delta_{c})}  \bigg]
\bigg\},
\end{equation*}
%--------------------------
else, repeat $\delta_{c}$. Again $f(\cdot)$ denotes the full likelihood
\eqref{eq:fulllikelihood} and $\pi(\beta,\sigma,\delta,\gamma)$ denotes the
priors.
%-------------------------------------------------
\item Sample $z_{i}|y, \beta, \nu_{i}, h_{i}, \sigma, \gamma, \delta \sim
    TN_{(\xi_{j-1},\,\xi_{j})}(x_{i}^{\prime}\beta + A\nu_{i}
    + C |\gamma | h_{i}, \sigma B\nu_{i})$ for $i =1, 2,\ldots,n$,
    where $\xi$ is obtained from $\delta$ by one-to-one mapping
    using equation~\eqref{eq:logtransformation}.
\end{enumerate}}
\rule{\textwidth}{0.5pt}
\end{algorithm}
\end{table*}
%-------------------------------------------------------------------------------

We next employ the Bayes' theorem and derive the joint posterior density as
proportional to the product of the likelihood and prior distributions. We
employ standard prior distributions as follows,
%-----------------------------------------
\begin{equation}
\begin{split}
\beta  & \sim N(\beta_{0}, B_{0}),  \hspace{0.5in}
\sigma \sim IG(n_{0}/2, d_{0}/2), \\
\gamma & \sim SB(L,U),                \hspace{0.65in}
\delta \sim N(\delta_{0}, D_{0}),
\end{split}
\label{eq:priors}
\end{equation}
%-----------------------------------------
where $N$, $IG$ and $SB$ denote normal, inverse-gamma and scaled-Beta
distributions, respectively. The lower and upper bounds of the scaled-Beta
distribution are obtained as mentioned in Section~\ref{sec:GAL}. Combining
the likelihood and the prior distributions, the augmented joint posterior
density can be written as,
%-----------------------------------------
\begin{equation}
\begin{split}
\pi(z,\beta,\nu,h,\sigma,\gamma,\delta|y) & \propto
f(y|z,\beta,\nu,h,\sigma,\gamma,\delta) \pi(z|\beta,\nu,h,\sigma,\gamma,\delta)
 \pi(\nu|\sigma)  \pi(h|\sigma) \pi(\beta)  \\
& \qquad \times \pi(\sigma)  \pi(\gamma) \pi(\delta) \\
%-----
& \propto \Big\{\prod_{i=1}^{n} f(y_{i}|z_{i},\beta,\nu_{i},h_{i},
\sigma,\gamma,\delta) \pi(\nu_{i}|\sigma) \pi(h_{i}|\sigma)\Big\}
\pi(z|\beta,\sigma,\nu,h,\gamma)  \\
& \qquad \times \pi(\beta) \pi(\sigma)  \pi(\gamma) \pi(\delta) \\
%-----
& \propto \Big\{\prod_{i=1}^{n} f(y_{i}|z_{i},\delta)
\pi(\nu_{i}|\sigma) \pi(h_{i}|\sigma)\Big\}
\pi(z|\beta,\sigma,\nu,h,\gamma)  \\
& \qquad \times \pi(\beta) \pi(\sigma)  \pi(\gamma) \pi(\delta),
\end{split}
\label{eq:jointPosterior}
\end{equation}
%-----------------------------------------
where the last line in the likelihood, based on $GAL(0,\sigma,\gamma)$, uses
the fact that given $z$ and $\delta$, the observed $y$ is independent of the
remaining parameters, because \eqref{eq:cutpoints} determines $y_{i}$ given
$(z,\delta)$ with probability 1. The conditional density of latent data $z$
is obtained from \eqref{eq:model3} and is given by
$\pi(z|\beta,\sigma,\nu,h,\gamma) = \prod_{i=1}^{n} N(z_{i}|x'_{i} \beta + A
\nu_{i}  + C |\gamma| h_{i}, \sigma B \nu_{i})$. Additionally, the prior
distributions for $(\beta,\sigma,\gamma,\delta)$ are assumed to be
independent in equation~\eqref{eq:jointPosterior}. Using the preceding
explanations, the ``complete data posterior'' in
equation~\eqref{eq:jointPosterior} can be expressed as,
%-----------------------------------------
\begin{equation}
\begin{split}
\pi(z,\beta,\nu,h,\sigma,\gamma,\delta|y) & \propto \bigg\{ \prod_{i=1}^{n}
1\big\{{\xi_{y_{i}-1}} < z_{i} \le \xi_{y_{i}} \big\}  N(z_{i}|x'_{i} \beta + A
\nu_{i}  + C |\gamma| h_{i}, \sigma B \nu_{i}) \\
%-----
& \qquad \times \mathcal{E}(\nu_{i}|\sigma) \, N^{+}(h_{i}|0,\sigma^{2})
\bigg\} N(\beta|\beta_{0},B_{0}) \, IG(\sigma|n_{0}/2,d_{0}/2) \\
%-----
& \qquad \times SB(\gamma|L,U) \, N(\delta|\delta_{0},D_{0}).
\end{split}
\label{eq:compDataPosterior}
\end{equation}
%-----------------------------------------
The objects of interest i.e., $(z,\beta,\nu,h,\sigma,\gamma,\delta)$ can be
sampled by deriving the conditional posterior densities from the complete
data posterior \eqref{eq:compDataPosterior} and judiciously using the full
likelihood \eqref{eq:fulllikelihood} as presented in
Algorithm~\ref{alg:algorithm1}. We note that our proposed algorithm is a form
of MH within partially collapsed Gibbs sampler and care has been taken to
guarantee convergence to stationary distribution as given in
\citet{vanDyk-Jiao-2015}.

Starting with the regression coefficients, $\beta$ is sampled from a normal
distribution, draws from which are programmed in most statistical software.
The scale and shape parameters $(\sigma,\gamma)$ are jointly sampled,
marginally of $(z,\nu,h)$, using a random-walk Metropolis-Hastings (MH)
algorithm with proposals drawn from a truncated bivariate normal
distribution. Joint sampling (together with the transformations $h_{i} =
\sigma s_{i}$ and $\nu_{i}=\sigma \omega_{i}$) is crucial for reducing the
high autocorrelation in MCMC draws as observed in \citet{Yan-Kottas-2017}.
The latent weight $\nu$ follows a generalized inverse-Gaussian (GIG)
distribution, draws from which are obtained using the technique proposed in
\citet{Devroye-2014}. Alternatively, one may employ the ratio of uniforms
method or the envelope rejection method
\citep{Dagpunar-1988,Dagpunar-1989,Dagpunar-2007}. The mixture variable $h$
is sampled from a half-normal distribution. Typical to ordinal models, the
cut-points $\delta$ do not have a tractable distribution and is sampled
marginally of $(z,\nu,h)$ using a random-walk MH algorithm
\citep[see][]{Jeliazkov-etal-2008,Rahman-2016}. Finally, the latent variable
$z$, conditional on the remaining parameters, is sampled from a truncated
normal distribution \citep{Botev-2017}. The derivations of the conditional
posteriors and details of the MH algorithms are presented in
\ref{app:FBQROR}.

%------------------------------------------------------------------------------
\section{Simulation Studies}\label{sec:simStudies}
%------------------------------------------------------------------------------
This section presents two simulation studies to demonstrate the performance
of the proposed algorithm and illustrate the advantages of the FBQROR model
compared to the BQROR model.

%-------------------------------------------------------------------------------
\subsection{Simulation Study~1}\label{sec:simStudy1}

In this simulation study, we estimate and compare the FBQROR model to the
BQROR model when errors are generated from a symmetric distribution.
Specifically, 300 observations are generated from the model $z_{i} = x'_{i}
\beta + \epsilon_{i}$, where covariates are sampled from a standard uniform
distribution $Unif[0,1]$, $\beta$ = $(2, -3, 4)'$ and $\epsilon$ is sampled
from a logistic distribution $\mathcal{L}(0,\pi^{2}/3)$. The resulting
continuous variable $z$ is symmetric and is utilized to construct the
discrete response variable $y$ based on the cut-point vector $\xi = (0, 2,
4)$. In our simulated data, the number of observations corresponding to the
four categories of $y$ are 42 (14\%) , 81 (27\%), 99 (33\%) and 78 (26\%),
respectively.

\begin{figure}[!t]
	\centerline{
		\mbox{\includegraphics[width=6.75in, height=4.5in]{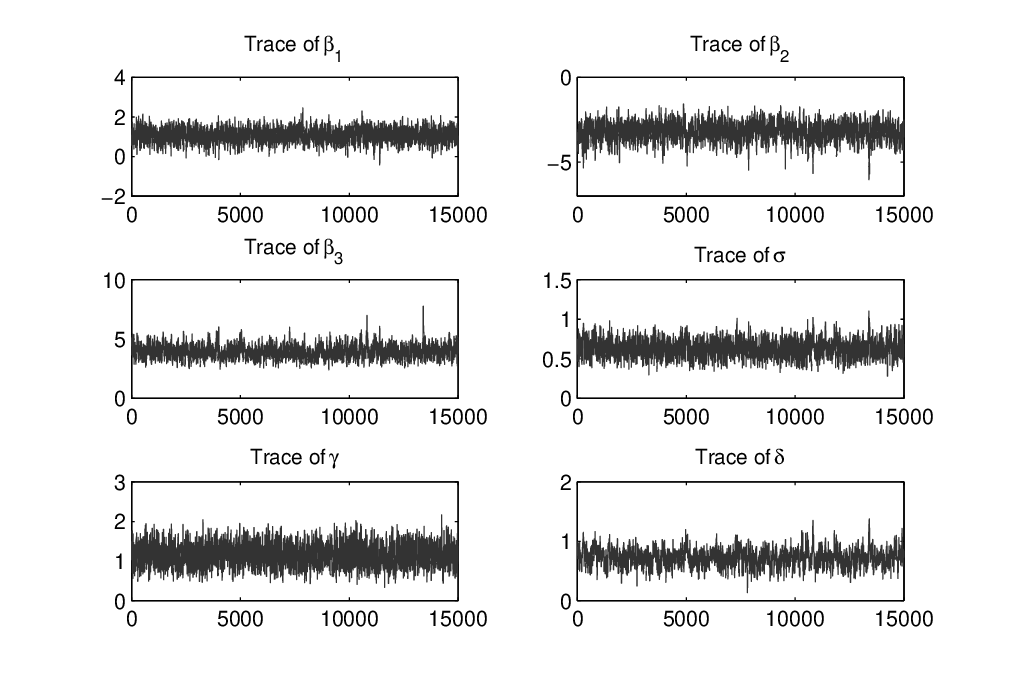}} }
	\caption{Trace plots of the MCMC draws at the 25th quantile for Simulation Study~1.}
\label{fig:SS1-MCMC-25th}
\end{figure}

The posterior estimates of the parameters for FBQROR model are obtained based
on the simulated data and the following moderately diffuse priors: $\beta
\sim N(0_{3}, 10 I_{3})$, $\sigma \sim IG(5/2,8/2)$, $\gamma \sim
SB(L,U,4,4)$ and $\delta \sim N(0_{J-3}, I_{J-3})$ for $p_{0}=(0.25, 0.5,
0.75)$, where $(L,U)$ depends on the value of $p_{0}$ as mentioned in
Section~\ref{sec:GAL}. Table~\ref{Table:SimResult1} reports the MCMC results
obtained from 15,000 iterations, after a burn-in of 5,000 iterations, along
with the inefficiency factors calculated using the batch-means method
\citep{Greenberg-2012}. The parameters $(\sigma, \gamma)$ are jointly sampled
using random-walk MH algorithm with tuning parameters $\iota_{1} =
(\sqrt{1.7}, \sqrt{2.25}, \sqrt{2.0})$ to achieve an acceptance rate of
approximately 33 percent for the three considered quantiles. Similarly,
$\delta$ is sampled using a random-walk MH algorithm with tuning factor
$\iota_{2}=(\sqrt{4.0}, \sqrt{3.2}, \sqrt{2.5})$ to obtain an acceptance rate
of around 33 percent. Inefficiency factors for all the model parameters are
low which imply low autocorrelation in MCMC draws and trace plots of the MCMC
iterations, as exhibited in Figure~\ref{fig:SS1-MCMC-25th} for the 25th
quantile, display quick convergence. Trace plots for the other two quantiles
are similar and have not been shown for the sake of brevity. The sampler is
reasonably quick and takes approximately 160 seconds per $1,000$ iterations.

%----------------------------  Table 1 ----------------------------------------
\begin{table}[!t]
\centering \footnotesize \setlength{\tabcolsep}{5pt} \setlength{\extrarowheight}{1.5pt}
\setlength\arrayrulewidth{1pt} \caption{Posterior mean (\textsc{mean}),
standard deviation (\textsc{std}) and inefficiency factor (\textsc{if}) of
the parameters in Simulation Study~1. The first panel presents results from
the FBQROR model and the second panel presents results from the BQROR model.}
\begin{tabular}{c rrr rrr rrr rrr r  }
\toprule
%------------------------------------------------------------------------------
& & \multicolumn{11}{c}{FBQROR Model} &  \\
\cmidrule{3-13}
& & \multicolumn{3}{c}{\textsc{25th quantile}} & & \multicolumn{3}{c}{\textsc{50th quantile}}
& & \multicolumn{3}{c}{\textsc{75th quantile}} \\
\cmidrule{3-5} \cmidrule{7-9}  \cmidrule{11-13}
$(\beta,\sigma,\gamma,\delta)$
                     & &  \textsc{mean} & \textsc{std} & \textsc{if}
                     & &  \textsc{mean} & \textsc{std} & \textsc{if}
                     & &  \textsc{mean} & \textsc{std} & \textsc{if} &  \\
\midrule
%-----------------------------------------------------------------------------
$\beta_{1}$     & & $ 1.07$  & $0.32$  & $2.94$  & & $ 2.17$ & $0.31$ & $2.66$
                & & $ 3.23$  & $0.35$  & $4.91$  \\
$\beta_{2}$     & & $-3.22$  & $0.50$  & $4.02$  & & $-3.14$ & $0.48$ & $3.67$
                & & $-3.13$  & $0.47$  & $4.72$  \\
%-----------------------------------------------------------------------------
$\beta_{3}$     & & $ 3.93$  & $0.53$  & $4.08$  & & $ 3.86$ & $0.51$ & $3.82$
                & & $ 3.90$  & $0.50$  & $5.28$  \\
%\rowcolor[gray]{0.92}
$\sigma$        & & $ 0.64$  & $0.10$  & $3.38$  & & $ 0.75$ & $0.09$ & $4.30$
                & & $ 0.60$  & $0.09$  & $4.10$  \\
%-----------------------------------------------------------------------------
$\gamma$        & & $ 1.14$  & $0.27$  & $2.55$  & & $-0.06$ & $0.17$ & $3.76$
                & & $-1.18$  & $0.24$  & $2.93$  \\
$\delta_{1}$    & & $ 0.73$  & $0.14$  & $4.75$  & & $ 0.71$ & $0.15$ & $4.84$
                & & $ 0.67$  & $0.14$  & $5.87$  \\
%-----------------------------------------------------------------------------
%------------------------------------------------------------------------------
\cmidrule{3-13}
& & \multicolumn{11}{c}{BQROR Model} &  \\
\cmidrule{3-13}
%------------------------------------------------------------------------------
& & \multicolumn{3}{c}{\textsc{25th quantile}} & & \multicolumn{3}{c}{\textsc{50th quantile}}
& & \multicolumn{3}{c}{\textsc{75th quantile}} &  \\
\cmidrule{3-5} \cmidrule{7-9}  \cmidrule{11-13}
%------------------------------------------------------------------------------
$(\beta,\sigma,\delta)$ & &  \textsc{mean} & \textsc{std} & \textsc{if}
                      & &  \textsc{mean} & \textsc{std} & \textsc{if}
                      & &  \textsc{mean} & \textsc{std} & \textsc{if} &   \\
\midrule
%------------------------------------------------------------------------------
$\beta_{1}$     & & $ 1.15$  & $0.30$  & $ 2.46$  & & $ 2.12$ & $0.29$ & $ 2.56$
                & & $ 2.86$  & $0.26$  & $ 4.43$  &  \\
$\beta_{2}$     & & $-3.23$  & $0.46$  & $ 3.21$  & & $-2.96$ & $0.43$ & $ 3.32$
                & & $-2.38$  & $0.46$  & $ 5.79$  &  \\
%------------------------------------------------------------------------------
$\beta_{3}$     & & $ 3.82$  & $0.44$  & $ 3.11$  & & $ 3.64$ & $0.47$ & $ 3.71$
                & & $ 3.01$  & $0.49$  & $ 6.46$  &  \\
$\sigma$        & & $ 0.65$  & $0.07$  & $ 3.55$  & & $ 0.71$ & $0.07$ & $ 3.84$
                & & $ 0.45$  & $0.05$  & $ 5.92$  &  \\
%------------------------------------------------------------------------------
$\delta_{1}$    & & $ 0.85$  & $0.13$  & $3.84$  & &  $ 0.60$ & $0.14$ & $ 4.51$
                & & $ 0.35$  & $0.15$  & $6.55$  \\
%------------------------------------------------------------------------------
\bottomrule
\end{tabular}
\label{Table:SimResult1}
\end{table}
%-------------------------------------------------------------------------------

%----------------------------  Figure 2 ----------------------------------------
%-------------------------------------------------------------------------------

The first panel of Table~\ref{Table:SimResult1} presents the results for the
FBQROR model. The results show that the posterior means for $\beta$ are close
to the true parameter values, the posterior mean of $\sigma$ adjusts the
scale of the distribution, and the posterior mean of $\delta_{1}$ yields a
value of $\xi_{3}$ close to 4, the true value used to generate the data. We
also estimate the BQROR model using a modification of Algorithm~1 in
\citet{Rahman-2016} -- by fixing the second cut-point and introducing a scale
parameter in the model. Prior distributions for $(\beta,\sigma,\delta)$ are
identical to that of the FBQROR model. The results, presented in the second
panel of Table~\ref{Table:SimResult1}, show that the posterior estimates for
$(\beta,\sigma,\delta)$ are similar. Moving to the skewness parameter in the
FBQROR model, the posterior mean of $\gamma$ at $p_{0}=(0.25, 0.50, 0.75)$
are $(1.14, -0.06,-1.18)$, which corresponds to a skewness of
$(0.01,0.20,0.04)$, respectively. Note that the posterior mean of $\gamma$ at
$p_{0}=0.50$ is statistically equivalent to zero. These skewness values imply
that the (latent) response variable is approximately symmetric at all
considered quantiles, which is reassuring since our data was generated from a
symmetric distribution. In contrast, the corresponding skewness values for
the BQROR model are $(1.64, 0, -1.64)$. Hence, the BQROR model fails to
accommodate the symmetric characteristic of the data at the 25th and 75th
quantiles.

%----------------------------  Table 2 ----------------------------------------
\begin{table}[!b]
\centering \footnotesize \setlength{\tabcolsep}{6pt}
\setlength{\extrarowheight}{1.5pt}
\setlength\arrayrulewidth{1pt} \caption{Model comparison using the conditional
log-likelihood ($\ln L$), Akaike information criterion (AIC)
and Bayesian information criterion (BIC) in Simulation Study~1.}
\begin{tabular}{l llll}
\toprule
&& \textsc{25th quantile}& \textsc{50th quantile} &\textsc{75th
quantile} \\
\cmidrule{3-3} \cmidrule{4-4} \cmidrule{5-5}
&&  \textsc{($\ln L$, AIC, BIC) }  &  \textsc{($\ln L$, AIC, BIC)}
 &  \textsc{($\ln L$, AIC, BIC) }
\\
\midrule		
%-----------------------------------------------------------------------------
\textsc{FBQROR}  & & $(-338, 688, 710)$  & $(-340, 691, 714)$
                   & $(-338, 688, 710)$  \vspace{0.2cm}  \\
\textsc{BQROR}   & & $(-346, 701, 720)$  & $(-340, 690, 708)$
                   & $(-348, 706, 724)$ \\		
%-----------------------------------------------------------------------------
\bottomrule
\end{tabular}
\label{Table:Sim1Modelfit}
\end{table}
%-----------------------------------------------------------------------------

We next investigate model fitness at different quantiles since various
choices of quantile $p_{0}$ may be interpreted as corresponding to a
different link function. Table~\ref{Table:Sim1Modelfit} presents the
conditional log-likelihood, the Akaike information criterion or AIC
\citep{Akaike-1974} and the Bayesian information criterion or BIC
\citep{Schwarz-1978} for both the FBQROR and BQROR models. Higher conditional
log-likelihood is preferable, while lower values of AIC/BIC indicate a better
model fit. As illustrated in Table~\ref{Table:Sim1Modelfit}, the conditional
log-likelihood for the FBQROR model is identical to the BQROR model at the
median, but higher at the other two considered quantiles. However, the FBQROR
model has an extra shape parameter and so to rule out the possibility of
higher log-likelihood arising due to additional parameters (i.e.,
overfitting), we compare the models using AIC and BIC. These two measures
introduce different penalty terms to account for the number of model
parameters. Based on AIC/BIC, there is strong evidence that the FBQROR model
provides a better fit at the 25th and 75th quantiles, but there is some
evidence in favor of BQROR model at the 50th quantile. The poor fit of the
BQROR model at the first and third quartiles reflects the rigidity of the AL
distribution, since $p_{0}=0.25$ ($0.75$) forces the AL distribution to be
positively (negatively) skewed.

%------------------------------------------------------------------------------
\subsection{Simulation Study 2}\label{sec:simStudy2}

%----------------------------  Table 3 ----------------------------------------
\begin{table}[!b]
\centering \footnotesize \setlength{\tabcolsep}{5pt} \setlength{\extrarowheight}{1.5pt}
\setlength\arrayrulewidth{1pt} \caption{Posterior mean (\textsc{mean}),
standard deviation (\textsc{std}) and inefficiency factor (\textsc{if}) of
the parameters in Simulation Study~2. The first panel presents results from
the FBQROR model and the second panel presents results from the BQROR model.}
\begin{tabular}{c rrr rrr rrr rrr r  }
\toprule
%------------------------------------------------------------------------------
& & \multicolumn{11}{c}{FBQROR Model} &  \\
%------------------------------------------------------------------------------
\cmidrule{3-13}
& & \multicolumn{3}{c}{\textsc{25th quantile}} & & \multicolumn{3}{c}{\textsc{50th quantile}}
& & \multicolumn{3}{c}{\textsc{75th quantile}} \\
\cmidrule{3-5} \cmidrule{7-9}  \cmidrule{11-13}
$(\beta,\sigma,\gamma,\delta)$  & &  \textsc{mean} & \textsc{std} & \textsc{if}
                     & &  \textsc{mean} & \textsc{std} & \textsc{if}
                     & &  \textsc{mean} & \textsc{std} & \textsc{if} &  \\
\midrule
%-----------------------------------------------------------------------------
$\beta_{1}$     & & $ 1.60$  & $0.36$  & $2.87$  & & $ 2.83$ & $0.37$ & $3.21$
                & & $ 4.50$  & $0.44$  & $4.00$  \\
$\beta_{2}$     & & $-6.50$  & $0.57$  & $3.37$  & & $-6.39$ & $0.61$ & $3.98$
                & & $-6.15$  & $0.64$  & $3.47$  \\
%-----------------------------------------------------------------------------
$\beta_{3}$     & & $ 3.69$  & $0.52$  & $3.01$  & & $ 3.59$ & $0.54$ & $3.36$
                & & $ 3.44$  & $0.58$  & $3.24$  \\
%\rowcolor[gray]{0.92}
$\sigma$        & & $ 0.74$  & $0.11$  & $8.08$  & & $ 0.74$ & $0.12$ & $2.25$
                & & $ 0.79$  & $0.10$  & $2.66$  \\
%-----------------------------------------------------------------------------
$\gamma$        & & $ 0.09$  & $0.15$  & $4.75$ & &  $-0.49$ & $0.09$ & $4.37$
                & & $-1.33$  & $0.17$  & $2.47$  \\
$\delta_{1}$    & & $ 0.91$  & $0.14$  & $2.20$  & & $ 0.86$ & $0.14$ & $2.89$
                & & $ 0.73$  & $0.13$  & $3.21$  \\
%------------------------------------------------------------------------------
\cmidrule{3-13}
& & \multicolumn{11}{c}{BQROR Model} &  \\
\cmidrule{3-13}
%------------------------------------------------------------------------------
& & \multicolumn{3}{c}{\textsc{25th quantile}} & & \multicolumn{3}{c}{\textsc{50th quantile}}
& & \multicolumn{3}{c}{\textsc{75th quantile}} &  \\
\cmidrule{3-5} \cmidrule{7-9}  \cmidrule{11-13}
%------------------------------------------------------------------------------
$(\beta,\sigma,\delta)$ & & \textsc{mean} & \textsc{std} & \textsc{if}
                      & &  \textsc{mean} & \textsc{std} & \textsc{if}
                      & &  \textsc{mean} & \textsc{std} & \textsc{if} &   \\
\midrule
%------------------------------------------------------------------------------
$\beta_{1}$     & & $ 1.10$  & $0.23$  & $ 2.72$  & & $ 1.91$ & $0.24$ & $ 2.44$
                & & $ 2.68$  & $0.24$  & $ 2.78$  &  \\
$\beta_{2}$     & & $-4.40$  & $0.40$  & $ 3.71$  & & $-3.88$ & $0.38$ & $ 2.84$
                & & $-3.26$  & $0.41$  & $ 3.20$  &  \\
%------------------------------------------------------------------------------
$\beta_{3}$     & & $ 2.49$  & $0.33$  & $ 3.01$  & & $ 2.11$ & $0.35$ & $ 2.75$
                & & $ 1.80$  & $0.34$  & $ 2.83$  &  \\
$\sigma$        & & $ 0.47$  & $0.04$  & $ 3.07$  & & $ 0.63$ & $0.06$ & $ 2.93$
                & & $ 0.49$  & $0.04$  & $ 2.94$  &  \\
%------------------------------------------------------------------------------
$\delta_{1}$    & & $ 0.51$  & $0.15$  & $3.23$  & &  $ 0.25$ & $0.14$ & $ 3.02$
                & & $-0.04$  & $0.14$  & $3.45$  \\
%------------------------------------------------------------------------------
\bottomrule
\end{tabular}
\label{Table:SimResult2}
\end{table}
%-------------------------------------------------------------------------------

Once again we estimate the FBQROR and BQROR models with simulated data, but
now the errors are generated from a chi-square distribution such that the
resulting distribution for the continuous latent variable $z$ is positively
skewed. In particular, 300 observations are generated from the model $z_{i} =
x'_{i} \beta + \epsilon_{i}$, where covariates are sampled from a standard
uniform distribution $Unif[0,1]$, $\beta$ = $(3, -7, 5)'$ and $\epsilon$ is
generated from $\chi^{2}(4) - 4$, i.e., a demeaned chi-square distribution.
The discrete response variable $y$ is obtained from $z$ based on cut-point
vector $\xi = (0, 3, 6)$, which yields 74 (24.67\%), 110 (36.67\%), 65
(21.67\%) and 51 (17.00\%) observations in the four categories of $y$.

\begin{figure}[!t]
	\centerline{
		\mbox{\includegraphics[width=6.750in, height=4.5in]{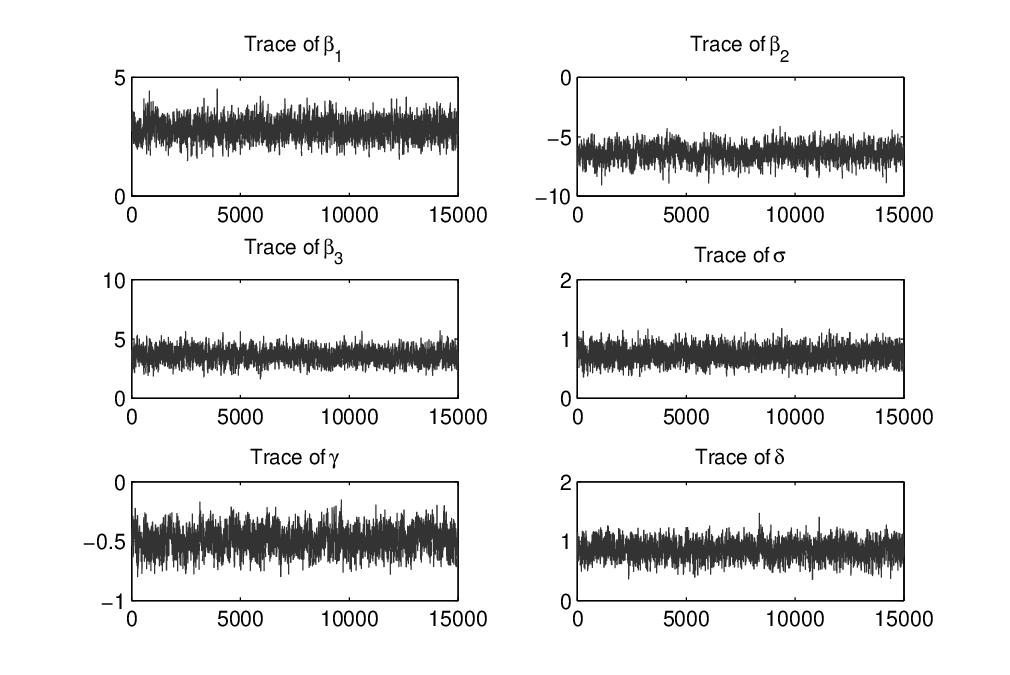}} }
	\caption{Trace plots of the MCMC draws at the 50th quantile for Simulation Study~2.}
\label{fig:SS2-MCMC-50th}
\end{figure}

Table~\ref{Table:SimResult2} reports the MCMC estimates obtained from 15,000
iterations after a burn-in of 5,000 iterations. Identical prior distributions
as in the first simulation study were used for both FBQROR and BQROR models.
The parameters $(\sigma, \gamma)$ and $\delta$ are sampled using random-walk
MH algorithm with tuning factors $\iota_{1} = (\sqrt{0.3},
\sqrt{0.7},\sqrt{1.45})$ and $\iota_{2}=(\sqrt{3.25},\sqrt{3.1},\sqrt{2.75})$
to achieve an acceptance rate of approximately 33 percent. The inefficiency
factors are low and trace plots, as displayed in
Figure~\ref{fig:SS2-MCMC-50th} for the 50th quantile, show quick convergence.
Trace plots at the other two quantiles are similar. Computational time
remains unchanged at approximately 160 seconds per $1,000$ iterations.

%----------------------------  Figure 3 ----------------------------------------
%-------------------------------------------------------------------------------

Table~\ref{Table:SimResult2} presents the results for the FBQROR and BQROR
models in the first and second panels, respectively. In the FBQROR model, the
posterior estimates of $\beta$ are close to the true values $(3,-7,5)$ and
the posterior estimates of $\sigma$ adjusts to capture the spread. Moreover,
the posterior estimates of $\gamma$ capture the skewness of the data
extremely well. Specifically, the posterior mean of $\gamma$ is not
statistically different from zero at the 25th quantile, so skewness is
approximately same as that of BQROR model (1.64). However, the skewness at
the 50th (75th) quantile is 1.31 ($0.16$) compared to a skewness of 0
($-1.64$) in the BQROR model. Hence, the FBQROR model correctly captures the
positive skewness of the simulated data at all quantiles. Model comparison
also points to the superiority of the FBQROR model as seen in
Table~\ref{Table:Sim2Modelfit}. The conditional log-likelihood for the FBQROR
model is higher than that of BQROR model at the 50th and 75th quantiles, but
identical at the 25th quantile. The AIC and BIC values suggest that there is
strong evidence to select the FBQROR model at the 50th and 75th quantiles,
and some evidence to favor the BQROR model at the 25th quantile. Once again,
the flexibility offered by the FBQROR model in terms of modeling the skewness
helps provide a better model fit compared to the rigid BQROR model.

%----------------------------  Table 4 ----------------------------------------
\begin{table}[!t]
\centering \footnotesize \setlength{\tabcolsep}{6pt}
\setlength{\extrarowheight}{1.5pt}
\setlength\arrayrulewidth{1pt} \caption{Model comparison using the conditional
log-likelihood ($\ln L$), Akaike information criterion (AIC)
and Bayesian information criterion (BIC) in Simulation Study~2.}
\begin{tabular}{l llll}
\toprule
&& \textsc{25th quantile}& \textsc{50th quantile} &\textsc{75th
quantile} \\
\cmidrule{3-3} \cmidrule{4-4} \cmidrule{5-5}
&&  \textsc{($\ln L$, AIC, BIC) }  &  \textsc{($\ln L$, AIC, BIC)}
 &  \textsc{($\ln L$, AIC, BIC) }
\\
\midrule		
%-----------------------------------------------------------------------------
\textsc{FBQROR}  & & $(-318, 649, 671)$  & $(-321, 654, 676)$
                   & $(-333, 678, 700)$  \vspace{0.2cm}  \\
\textsc{BQROR}   & & $(-318, 646, 664)$  & $(-331, 673, 692)$
                   & $(-358, 727, 745)$ \\		
%-----------------------------------------------------------------------------
\bottomrule
\end{tabular}
\label{Table:Sim2Modelfit}
\end{table}
%-----------------------------------------------------------------------------

%------------------------------------------------------------------------------
\section{Application}\label{sec:application}
%------------------------------------------------------------------------------

In the US, consumers have typically viewed homeownership as a good long-term
financial investment. However, the recent housing crisis and the subsequent
economic recession (Dec 2007 - Jun 2009) had substantial adverse effects,
particularly on homeowners. Housing values depreciated considerably, more
than four million foreclosures took place between 2008-2011, and one-fourth
of all homeowners were living in houses worth less than the mortgage at the
peak of the crisis \citep{Belsky-2013}. Homeownership rate declined from 69.2
percent during Q2, 2004 to 66.4 percent during Q1, 2011 and further to 62.9
percent during Q2, 2016 (Source: US Bureau of the Census). These adverse
effects may have fundamentally altered the perceived benefits of
homeownership as a good long-term investment (see \citet{Rohe-Lindblad-2014}
for a conceptual model). Consequently, it is of considerable interest to
analyze public opinion on homeownership as an investment and examine how
socioeconomic factors, demographic variables and exposure to financial
distress affect public responses.

%----------------------------  Table 5 ----------------------------------------
\begin{table}[!b]
\centering \footnotesize \setlength{\tabcolsep}{6pt} \setlength{\extrarowheight}{1.5pt}
\setlength\arrayrulewidth{1pt}
\caption{Descriptive summary of the variables.}
\begin{tabular}{lp{8.5cm}r r}
\toprule
\textsc{variable}      & \textsc{description}   & \textsc{mean}   & \textsc{std}  \\
\midrule
%------------------------------------------------------------------------------
\textsc{log age}       & Logarithm of age (in years)      &  3.71   &  0.44       \\
\textsc{log income}    & Logarithm of the mid-point of income category (in US dollars)
                                                          & 10.68   &  0.95       \\
\textsc{household size}& Number of members in the household
                                                          &  2.92   &  1.66       \\
\midrule
%------------------------------------------------------------------------------
                        &     & \textsc{count} & \textsc{percent}  \\
                        \midrule
\textsc{female}         & Indicator variable for female gender
                                                          &    925   &    51.42   \\
\textsc{post-bachelors} & Respondent's highest qualification is Masters,
                          Professional or Doctorate       &    257   &    14.29   \\
\textsc{bachelors}      & Respondent's highest qualification is Bachelors
                                                          &    395   &    21.96   \\
\textsc{below bachelors}& Respondent holds a 2-year associate degree, went to
                          some college with no degree, or attended technical,
                          trade or vocational school after high school
                                                          &    551   &    30.63   \\
\textsc{hs and below}   & Respondent is a high school (HS)
                          graduate or below               &    596   &    33.92   \\
%------------------------------------------------------------------------------
\textsc{full-time}      & Works full time                 &    849   &    47.19   \\
\textsc{part-time}      & Works part time                 &    266   &    14.79   \\
\textsc{Unemployed}     & Either unemployed, student or retired
                                                          &    684   &    38.02   \\
%------------------------------------------------------------------------------
\textsc{white}          & Respondent is a White-American
                                                          &   1293   &    71.87   \\
\textsc{african-american}
                        & Respondent is an African-American
                                                          &    272   &    15.12   \\
\textsc{all other races}
                        & Respondent is an Asian, Asian-American or belongs to
                          some other race
                                                          &    234   &    13.01   \\
%------------------------------------------------------------------------------
\textsc{northeast}      & Lives in the northeast region of US
                                                          &    249   &    13.84   \\
\textsc{west}           & Lives in the west region of US
                                                          &    408   &    22.68   \\
\textsc{south}          & Lives in the south region of US
                                                          &    822   &    45.69   \\
\textsc{midwest}        & Lives in the midwest region of US
                                                          &    320   &    17.79   \\
%------------------------------------------------------------------------------
\textsc{fin-better}     & Financially better-off post the Great Recession
                                                          &    537   &    29.85   \\
\textsc{fin-same}       & Financially equivalent pre and post the Great Recession
                                                          &    445   &    24.74   \\
\textsc{fin-worse}      & Financially worse-off post the Great Recession
                                                          &    817   &    45.41   \\
%------------------------------------------------------------------------------
\midrule
                        & Strongly disagree or somewhat disagree that
                          homeownership is the best long-term investment (LTI)
                          in US                           &   310    &   17.23     \\
\textsc{opinion}        & Somewhat agree that homeownership is the best LTI
                                                          &   828    &   46.03    \\
                        & Strongly agree that homeownership is the best LTI
                                                          &   661    &   36.74    \\
%------------------------------------------------------------------------------
\bottomrule
\end{tabular}
\label{Table:DataSummary}
\end{table}
%------------------------------------------------------------------------------

This paper utilizes the Higher Education/Housing Survey data of March 2011,
conducted by the Princeton Survey Research Associates International and
sponsored by the Pew Social and Demographic Trends project. Interviews were
conducted over the telephone between March 15-29, 2011 on a nationally
representative sample of 2,142 adults living in the continental US. After
removing missing responses, we are left with a sample of 1,799 observations
for our analysis. Our dependent variable is response to the statement, ``Some
people say that buying home is the best long-term investment in the United
States. Do you strongly agree, somewhat agree, somewhat disagree or strongly
disagree?'' Responses are recorded into one of the four categories, however,
we append the responses ``strongly disagree'' and ``somewhat disagree'' as
the former category had less than 5 percent observations. The survey also
collected information on a wide range of socioeconomic, demographic and
geographic variables, some of which are used as covariates in the model.
Table~\ref{Table:DataSummary} presents the description and summary statistics
of all covariates and the response variable utilized in the study.

The average age of the sampled individuals is 44.84 years with a standard
deviation of 18.59 years. Information on family income in the survey is
recorded as one of 9 income categories: $ < 10k$, $10k-20k$, $20k-30k$,
$30k-40k$, $40k-50k$, $50k-75k$, $75k-100k$, $100k-150k$ and $>150k$, where
$k$ denotes a thousand dollars and \$5,000 and \$1,70,000 have been imputed
for the first and last income categories. We include the logarithm of the
mid-point of the income category as a variable in the model. Mean household
size is 2.92 with a standard deviation of 1.66. The percent of female is
slightly more than that of males. Educational classification shows that HS
and below forms the largest category (33.92\%) and post-bachelors forms the
smallest category (14.29\%), with proportions decreasing as we move from the
lowest to the highest educational category. Employment status shows that
61.98\% are either employed full-time or part-time, while the remaining are
unemployed, student or retired. With respect to race, the sample is
predominantly white (71.87\%), followed by African-Americans (15.12\%) and
all other races (13.01\%). Geographically, most of the sampled individuals
live in the South (45.69\%), followed by West (22.68\%), Midwest (17.79\%)
and Northeast (13.84\%). These regional classifications are as defined by the
US Census Bureau. To measure exposure to financial distress, we include self
reported financial condition pre and post Great Recession. As expected,
almost half the sampled individuals (45.41\%) are financially worse-off post
the Great Recession.

Moving to the response variable, Table~\ref{Table:DataSummary} shows that
more than three-fourths of the sampled individuals ($82.77\%$) somewhat or
strongly agree that homeownership is the best long-term investment.
Therefore, US public opinion on homeownership remains largely unchanged even
after the housing meltdown and the Great Recession. A similar conclusion has
been obtained using data from the Survey of Consumers collected by the
University of Michigan and the National Housing Survey collected by Fannie
Mae \citep{Belsky-2013}. This result is primarily due to the financial
benefits of homeownership making owning more lucrative than renting,
especially in the long run. Two related articles that have studied the
preference for homeownership versus renting using binary models on survey
data are \citet{Bracha-Jamison-2012} and \citet{Drew-Herbert-2013}. Both
studies find no fundamental shifts in attitude towards homeownership.

\begin{figure}[!t]
	\centerline{
		\mbox{\includegraphics[width=6.75in, height=8.0in]{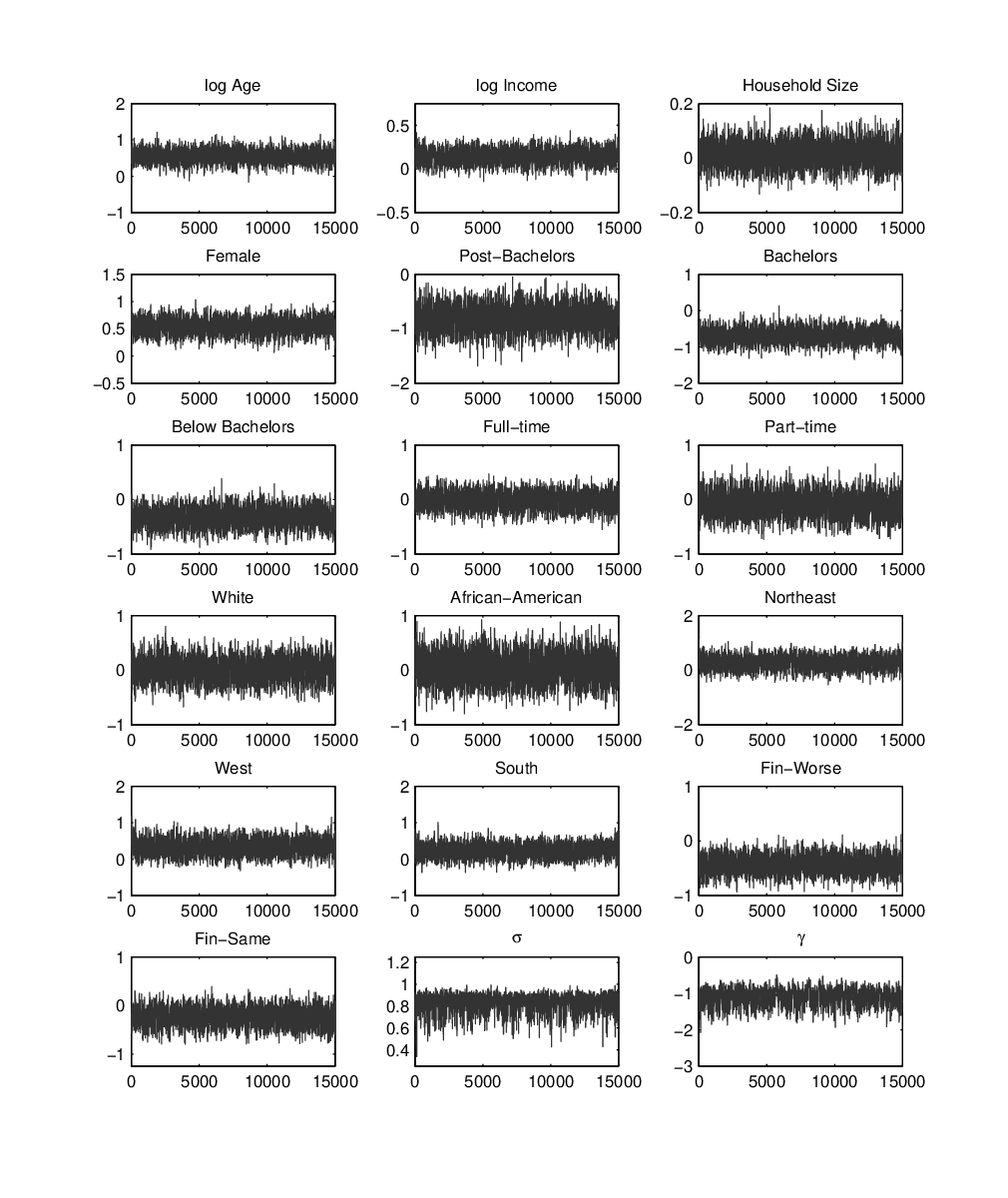}} }
\vspace{-2pc}
	\caption{Trace plots of the MCMC draws at the 75th quantile for the homeownership
             application.}
\label{fig:HomeO-MCMC-75th}
\end{figure}

We employ the FBQROR and BQROR models to analyze public opinion on
homeownership as the best long-term investment based on the covariates
presented in Table~\ref{Table:DataSummary}. The MCMC results, presented in
Table~\ref{Table:HomeOwnResults}, are based on 15,000 iterations after a
burn-in of 5,000 iterations with identical priors as in the simulation
studies. With three values of the ordinal response variable, we have two
cut-points and they are fixed at $(0,3)$ for both models across quantiles.
Similar to the simulation studies, $(\sigma,\gamma)$ is sampled using joint
random-walk MH algorithm with tuning parameters $\iota_{1}=(\sqrt{3.0},
\sqrt{2.4}, \sqrt{4.4})$ to get an acceptance rate of approximately 33
percent for the three considered quantiles. The inefficiency factor of the
parameters are all less than 5 and trace plots of MCMC draws, as displayed in
Figure~\ref{fig:HomeO-MCMC-75th} for the 75th quantile, show good mixing.
Trace plots at the other two quantiles are similar.

%----------------------------  Table 6 ----------------------------------------
\begin{table}[!t]
\centering \footnotesize \setlength{\tabcolsep}{2.5pt} \setlength{\extrarowheight}{1.75pt}
\setlength\arrayrulewidth{1pt}
\caption{Posterior mean (\textsc{mean}) and standard deviation
(\textsc{std}) of the parameters in the FBQROR and
BQROR models for the homeownership application.}
\begin{tabular}{l rrr rrr rrr rrr rrrr rr  }
\toprule
& & \multicolumn{8}{c}{\textsc{fbqror}} && \multicolumn{8}{c}{\textsc{bqror}}  \\
\cmidrule{3-10}  \cmidrule{12-19}
%------------------------------------------------------------------------------
&& \multicolumn{2}{c}{\textsc{25th}} && \multicolumn{2}{c}{\textsc{50th}}
&& \multicolumn{2}{c}{\textsc{75th}} && \multicolumn{2}{c}{\textsc{25th}}
&& \multicolumn{2}{c}{\textsc{50th}} && \multicolumn{2}{c}{\textsc{75th}}   \\
\cmidrule{3-4} \cmidrule{6-7}  \cmidrule{9-10} \cmidrule{12-13}   \cmidrule{15-16} \cmidrule{18-19}
%------------------------------------------------------------------------------
                   &&  \textsc{mean} & \textsc{std} &&  \textsc{mean} & \textsc{std}
                   &&  \textsc{mean} & \textsc{std} &&  \textsc{mean} & \textsc{std}
                   &&  \textsc{mean} & \textsc{std} &&  \textsc{mean} & \textsc{std}     \\
\midrule
\textsc{intercept}       && $ -3.11$  & $0.93$  && $ -1.72$ & $1.02$  && $  0.18$ & $0.90$
                         && $ -3.08$  & $0.89$  && $ -1.43$ & $0.93$  && $  1.44$ & $0.62$  \\
%--------------------------
\textsc{log age}         && $  0.52$  & $0.17$  && $  0.60$ & $0.17$  && $  0.55$ & $0.16$
                         && $  0.18$  & $0.16$  && $  0.63$ & $0.17$  && $  0.47$ & $0.12$  \\
\textsc{log income}      && $  0.18$  & $0.08$  && $  0.16$ & $0.09$  && $  0.15$ & $0.08$
                         && $  0.25$  & $0.08$  && $  0.14$ & $0.08$  && $  0.04$ & $0.05$  \\
%--------------------------
\textsc{household size}  && $  0.01$  & $0.04$  && $  0.01$ & $0.04$  && $  0.01$ & $0.04$
                         && $ -0.01$  & $0.04$  && $  0.01$ & $0.04$  && $  0.02$ & $0.03$  \\
\textsc{female}          && $  0.64$  & $0.14$  && $  0.62$ & $0.14$  && $  0.54$ & $0.13$
                         && $  0.57$  & $0.13$  && $  0.57$ & $0.13$  && $  0.32$ & $0.09$  \\
%--------------------------
\textsc{post-bachelors}  && $ -0.83$  & $0.22$  && $ -0.86$ & $0.22$  && $ -0.81$ & $0.20$
                         && $ -0.50$  & $0.21$  && $ -0.86$ & $0.21$  && $ -0.62$ & $0.15$  \\
\textsc{bachelors}       && $ -0.72$  & $0.20$  && $ -0.74$ & $0.19$  && $ -0.69$ & $0.18$
                         && $ -0.43$  & $0.18$  && $ -0.74$ & $0.18$  && $ -0.52$ & $0.13$  \\
%--------------------------
\textsc{below bachelors} && $ -0.37$  & $0.17$  && $ -0.36$ & $0.17$  && $ -0.33$ & $0.16$
                         && $ -0.26$  & $0.16$  && $ -0.34$ & $0.16$  && $ -0.25$ & $0.12$  \\
\textsc{full-time}       && $ -0.04$  & $0.16$  && $ -0.02$ & $0.16$  && $ -0.02$ & $0.14$
                         && $ -0.09$  & $0.15$  && $ -0.01$ & $0.15$  && $  0.01$ & $0.10$  \\
%--------------------------
\textsc{part-time}       && $ -0.07$  & $0.21$  && $ -0.05$ & $0.21$  && $ -0.09$ & $0.19$
                         && $ -0.01$  & $0.19$  && $ -0.06$ & $0.20$  && $ -0.12$ & $0.14$  \\
\textsc{white}           && $  0.04$  & $0.21$  && $  0.01$ & $0.21$  && $  0.01$ & $0.18$
                         && $  0.18$  & $0.18$  && $ -0.03$ & $0.20$  && $ -0.06$ & $0.14$  \\
%--------------------------
\textsc{african-american}&& $  0.05$  & $0.26$  && $  0.03$ & $0.27$  && $  0.05$ & $0.24$
                         && $  0.17$  & $0.23$  && $ -0.01$ & $0.25$  && $ -0.01$ & $0.17$  \\
\textsc{northeast}       && $  0.37$  & $0.24$  && $  0.34$ & $0.24$  && $  0.27$ & $0.22$
                         && $  0.41$  & $0.22$  && $  0.28$ & $0.23$  && $  0.12$ & $0.16$  \\
%--------------------------
\textsc{west}            && $  0.44$  & $0.22$  && $  0.46$ & $0.22$  && $  0.34$ & $0.19$
                         && $  0.39$  & $0.19$  && $  0.42$ & $0.21$  && $  0.17$ & $0.14$  \\
\textsc{south}           && $  0.34$  & $0.19$  && $  0.35$ & $0.19$  && $  0.23$ & $0.17$
                         && $  0.38$  & $0.17$  && $  0.30$ & $0.18$  && $  0.05$ & $0.12$  \\
%--------------------------
\textsc{fin-worse}       && $ -0.45$  & $0.16$  && $ -0.48$ & $0.16$  && $ -0.43$ & $0.14$
                         && $ -0.34$  & $0.15$  && $ -0.46$ & $0.15$  && $ -0.31$ & $0.11$  \\
\textsc{fin-same}        && $ -0.25$  & $0.18$  && $ -0.26$ & $0.18$  && $ -0.24$ & $0.17$
                         && $ -0.25$  & $0.18$  && $ -0.25$ & $0.17$  && $ -0.15$ & $0.13$  \\
%--------------------------
$\sigma$                 && $  0.89$  & $0.13$  && $  1.07$ & $0.05$  && $  0.83$ & $0.08$
                         && $  0.88$  & $0.03$  && $  1.06$ & $0.04$  && $  0.59$ & $0.02$  \\
$\gamma$                 && $  1.05$  & $0.31$  && $ -0.10$ & $0.10$  && $ -1.09$ & $0.24$
                         && $..$      & $..$    && $..$     & $..$    && $..$     & $..$    \\
%------------------------------------------------------------------------------------------------
\bottomrule
\end{tabular}
\label{Table:HomeOwnResults}
\end{table}
%------------------------------------------------------------------------------

%----------------------------  Figure 4 ----------------------------------------
%\afterpage{\FloatBarrier}
%-------------------------------------------------------------------------------

The results presented in Table~\ref{Table:HomeOwnResults} clearly show that
the posterior estimates from FBQROR and BQROR models are fairly similar
across all quantiles. Hence, we restrict our attention to the FBQROR model
and use the BQROR model for model comparison. Moreover, we primarily discuss
the covariates which are statistically different from zero at the 95\%
probability level. As seen from Table~\ref{Table:HomeOwnResults}, age has a
positive effect which implies that older individuals are more likely to
strongly agree that homeownership is the best long-term investment. This
result is consistent with the view that older adults are less likely to
change their attitude when faced with harsh economic experiences such as an
economic crisis \citep{Malmendier-2011,Giuliano-Spilimbergo-2014}. Our result
also finds support in \citet{Bracha-Jamison-2012}, who find that older
individuals are more confident about homeownership (relative to renting)
following large price declines. Income has a positive effect implying that
higher income individuals are more likely to agree with the investment
benefits of homeownership. However, income is an important factor only at the
25th quantile. This result is somewhat consistent with
\citet{Drew-Herbert-2013}, who find no statistically significant association
between income and viewing homeownership as a better financial choice over
renting.

Opinions across gender often vary due to risk perceptions and this is
reflected in our results. We find that females are more likely to strongly
agree that homeownership is the best long-term investment across all
quantiles. This result is consistent with the view that females are more risk
averse than males and homeownership has historically been a safe investment.
However, our results are in contrast to \citet{Bracha-Jamison-2012}, who find
that females are more uncertain about the financial gain from buying a house.
Higher education has a negative effect on positive opinion about
homeownership. The negative post-bachelors coefficient indicates that an
individual with a post-bachelors degree (relative to HS and below education)
will be less willing to strongly agree that homeownership is the best
long-term investment. Similarly, individuals with a bachelors or below
bachelors education are less likely to positively view the investment
benefits of homeownership. The negative effect of higher education on home
ownership is also reported in \citet{Bracha-Jamison-2012} and
\citet{Drew-Herbert-2013}.

Employment status, whether full-time or part-time as compared to being
unemployed, does not account for differences in opinion on the financial
benefits of homeownership. The same is true for the White and
African-American race indicators. Thus, individuals have similar views on
homeownership as an investment irrespective of employment status or race.
During the housing crisis, the decline in housing prices varied tremendously
across geographic regions. West and South regions experienced the largest
declines in housing prices. Accordingly, we include indicator variables for
geographic regions to capture differences in opinion due to residing  in
different regions. The results from the lower quantiles suggest that
individuals living in the West, relative to Midwest, are more likely to
strongly agree on the financial benefits of homeownership. This result is
interesting since people living in the West were the hardest hit in terms of
housing price declines. Undoubtedly, the housing meltdown and the economic
crisis caused serious financial distress to a large number of individuals in
the US. This experience may have altered views on homeownership. To capture
the effect of financial distress on homeownership views, we include indicator
variables for post recession financial situation. The results indicate that
individuals who are financially worse-off post the Great Recession, relative
to those who are better-off, are less probable to strongly agree that
homeownership is the best long-term investment. Hence, our results provide
evidence that financial hardship endured during the Great Recession
negatively impacted public views on the investment benefits of homeownership.

%----------------------------  Table 7 ----------------------------------------
\begin{table}[!t]\centering
\footnotesize \setlength{\tabcolsep}{3.5pt} \setlength{\extrarowheight}{2pt}
\caption{Change in predicted probabilities of the responses: somewhat or
strongly disagree $(y=1)$, somewhat agree $(y=2)$ and strongly agree $(y=3)$ that
homeownership is the best long-term investment.}
\begin{tabular}{lrr rrr rrr rrr r}
\toprule
    & & \multicolumn{3}{c}{\textsc{female}}
    & & \multicolumn{3}{c}{\textsc{post-bachelors}}
    & & \multicolumn{3}{c}{\textsc{finc-worse}} \\
    \cmidrule{3-5} \cmidrule{7-9}  \cmidrule{11-13}
%-----------------------------------------------------------------------------------------
    & &  25th   &  50th   & 75th    & &  25th     &  50th    &  75th
    & &  25th   &  50th   & 75th         \\ [0.5ex] \hline \\[-6pt]
$\Delta$P(y=1)
    &&  $-0.0630$  &  $-0.0553$   & $-0.0504$    &&  $ 0.0962$  &  $ 0.0900$   & $ 0.0852$
    &&  $ 0.0452$  &  $ 0.0429$   & $ 0.0404$    \\
$\Delta$P(y=2)
    &&  $-0.0322$  &  $-0.0405$   & $-0.0467$    &&  $ 0.0174$  &  $ 0.0282$   & $ 0.0473$
    &&  $ 0.0219$  &  $ 0.0302$   & $ 0.0360$    \\
$\Delta$P(y=3)
    &&  $ 0.0952$  &  $ 0.0958$   & $ 0.0970$    &&  $-0.1136$  &  $-0.1182$   & $-0.1325$
    &&  $-0.0671$  &  $-0.0731$   & $-0.0764$    \\
\bottomrule
\end{tabular}
\label{Table:HOCovEffect}
\end{table}
%------------------------------------------------------------------------------------------

In the previous paragraphs, we discussed the direction of covariate effects
on the last outcome i.e., strongly agree that homeownership is the best
long-term investment. The direction of covariate effect on the first outcome
(strongly disagree or somewhat disagree) is the opposite, while the effect on
the second outcome (somewhat agree) cannot be known \emph{a-priori}. This is
because the link function in ordinal models is non-linear and hence the
regression coefficients do not give the covariate effects. To make it clear,
we calculate the marginal effect for three variables: female, post-bachelors
and worse financial condition \citep{Jeliazkov-Vossmeyer-2018}. The change in
predicted probabilities for the three response are reported in
Table~\ref{Table:HOCovEffect}. We see that at the 25th quantile, individuals
who are exposed to financial distress (i.e., financially worse-off) are
6.71\% less likely to `strongly agree', 2.19\% more likely to `somewhat
agree' and 4.52\% more likely to `strongly disagree or somewhat disagree'
that homeownership is the best long-term investment. The marginal effect of
financial distress on the responses `strongly agree' are more pronounced at
the 50th and 75th quantiles. We can similarly interpret the change in
predicted probabilities for female and post-bachelors education on the three
responses for different quantiles.

%----------------------------  Table 8 ----------------------------------------
\begin{table}[!b]
\centering \footnotesize \setlength{\tabcolsep}{6pt}
\setlength{\extrarowheight}{1.5pt}
\setlength\arrayrulewidth{1pt} \caption{Model comparison using the conditional
log-likelihood ($\ln L$), Akaike information criterion (AIC)
and Bayesian information criterion (BIC) in the homeownership application.}
\begin{tabular}{l llll}
\toprule
&& \textsc{25th quantile}& \textsc{50th quantile} &\textsc{75th
quantile} \\
\cmidrule{3-3} \cmidrule{4-4} \cmidrule{5-5}
&&  \textsc{($\ln L$, AIC, BIC) }  &  \textsc{($\ln L$, AIC, BIC)}
 &  \textsc{($\ln L$, AIC, BIC) }
\\
\midrule		
%-----------------------------------------------------------------------------
\textsc{FBQROR}  & & $(-1816, 3671, 3775)$  & $(-1815, 3668, 3772)$
                   & $(-1816, 3671, 3775)$  \vspace{0.2cm}  \\
\textsc{BQROR}   & & $(-1824, 3684, 3783)$  & $(-1815, 3665, 3764)$
                   & $(-1818, 3673, 3772)$ \\		
%-----------------------------------------------------------------------------
\bottomrule
\end{tabular}
\label{Table:HomeOModelfit}
\end{table}
%-----------------------------------------------------------------------------

To assess model fitness across quantiles, we report the conditional
log-likelihood, AIC and BIC in Table~\ref{Table:HomeOModelfit}. The
log-likelihood for the FBQROR model across all quantiles are higher or
similar to that obtained from the BQROR model. However, according to AIC
there is strong (weak) evidence in favor of the FBQROR model at the 25th
(75th) quantile, but weak evidence in favor of the BQROR model at the 50th
quantile. Based on BIC, there is strong evidence to prefer the FBQROR (BQROR)
model at the 25th (50th) quantile, but positive evidence to prefer the BQROR
model at the 75th quantile.

%------------------------------------------------------------------------------
\section{Conclusion}\label{sec:conclusion}
%------------------------------------------------------------------------------

This paper presents an estimation algorithm for Bayesian quantile regression
in univariate ordinal models where the error is assumed to follow a GAL
distribution, referred to as the FBQROR model. To propose this estimation
procedure, we explore the GAL distribution and introduce and derive its
cumulative distribution function and moment generating function. We show that
the advantages offered by the GAL distribution --- which allows the mode,
skewness and tails to vary for any given quantile --- can be gainfully
utilized to better estimate Bayesian quantile regression in ordinal models.
We also emphasize on the efficiency of the MCMC algorithm, which is attained
through suitable transformation of the variables and joint sampling of the
scale and shape parameters. The practical advantages of the proposed model
are illustrated in multiple simulation studies via model comparison, where it
is observed that the FBQROR model can provide a better model fit compared to
an ordinal model with an AL distribution, labeled BQROR model
\citep{Rahman-2016}. Our proposed algorithm is also implemented to examine US
public opinion on homeownership as the best long-term investment following
the Great Recession. The results provide interesting insights which may be
useful for policymakers and researchers on US housing market.

The GAL distribution proposed in \citet{Yan-Kottas-2017} and further studied
in this paper is relatively new and needs to be studied further, particularly
due to its usefulness in Bayesian quantile regression. In fact, the GAL
distribution can practically be employed to estimate most Bayesian quantile
regression models that have been estimated using the AL distribution. A
partial list includes the Tobit model with endogenous covariates, censored
model, count data model, mixed-effect or longitudinal data model (work in
progress) and censored dynamic panel data model. Moreover, the distribution
can also be utilized to explore Bayesian variable selection in all the above
mentioned models. We leave these opportunities for future research.

%------------------------------------------------------------------------------
%-------------------- Beginning of Appendices ---------------------------------
%------------------------------------------------------------------------------
\clearpage
\newpage \appendix \renewcommand\thesection{Appendix \Alph{section}}
\setcounter{equation}{0}
\renewcommand\theequation{A.\arabic{equation}}
%\setcounter{secnumdepth}{-1} % This code removes the numbering of the appendix
%-------------------------------------------------------------------------------
\section{The GAL Distribution} \label{app:GAL}

This appendix derives the \emph{pdf} of the GAL distribution from the mixture
representation of the AL distribution, and introduces the \emph{cdf} and
\emph{mgf} of the GAL distribution. The \emph{mgf} is also utilized to derive
the mean, variance and skewness of the GAL distribution.

%------------------------------------------------------------------------------
\subsection{Probability Density Function}\label{app:GALpdf}

\textbf{Theorem 1}: Suppose $Y \sim GAL(\mu, \sigma, p, \alpha)$ and has the
\emph{pdf} given by equation~\eqref{eq:GALpdf}, then Y has the following
hierarchical representation, $Y = \mu + \alpha \sigma S + \sigma A(p) W +
\sigma [B(p) W]^{\frac{1}{2}} U$, where all the notations are as in
Section~\ref{sec:GAL}.

\textbf{Proof}: Using the mixture representation we can write the \emph{pdf}
of $Y$ as,
\begin{align}
f(y|\theta) & =
\int_{\mathbb{R^+}}\int_{\mathbb{R^+}} N(y | \mu + \sigma\alpha s + \sigma
A(p)w, \sigma^2B(p)w )\exp(w|1) N^{+}(s | 0,1) dw \, ds \nonumber\\
%---------------------------
& =  \int_{\mathbb{R^+}} \frac{1}{\sqrt{\sigma^2B(p)}} \;
     \int_0^\infty \left[ \frac{1}{\sqrt{2\pi w}}\exp{\left \{-\frac{1}{2}
     \left[ \frac{(y-b-aw)^2}{cw} \right] -w \right \} dw }\right]
     N^+(s | 0,1) \; ds \label{eq:pdfDblInt}
\end{align}
%---------------------------
where $\bm{\theta} = (\mu,\sigma,p,\alpha)$, $a=\sigma
A(p)=\frac{\sigma(1-2p)}{p(1-p)}$, $b=\mu +\sigma\alpha s$, and
$c=\sigma^2B(p)=\frac{2\sigma^2}{p(1-p)}$. We let $P$ denote the second
integral and integrate with respect to $w$ as follows,
%---------------------------
\begin{align}
P & = \int_0^\infty \frac{1}{\sqrt{2\pi w}} \exp\left\{ -\frac{1}{2}
\left[\frac{(a^2+2c)w^2+(y-b)^2-2a(y-b) w)}{c w} \right] \right\}  dw
\nonumber \\
& =
\exp{\left\{ \frac{a(y-b)}{c} \right\} } \int_0^\infty \frac{1}{\sqrt{2\pi w}}
\exp\left\{ -\frac{1}{2} \left[ \frac{(a^2+2c)}{c}w+\frac{(y-b)^2}{cw}
\right]\right\} dw \nonumber \\
& =
\exp{\left\{ \frac{a(y-b)}{c} \right\}}\int_0^\infty \sqrt{\frac{\gamma}{2\pi w}}
\, \sqrt{\frac{1}{\gamma}} \exp\left\{ -\frac{1}{2} \left[ \gamma w +
\frac{\gamma}{\mu^{2} w} - \frac{2\gamma}{\mu} + \frac{2\gamma}{\mu}
\right] \right\} dw \nonumber \\
& =
\frac{1}{\sqrt{\gamma}} \exp{\left\{ -\frac{\gamma}{\mu} \right\} }\exp{\left\{
\frac{a(y-b)}{c} \right\}} \underbrace{\int_0^\infty \sqrt{\frac{\gamma}{2\pi w}}
\exp \left\{ -\frac{\gamma}{2} \left( \frac{(1-\mu w)^2}{\mu^2 w}\right)
\right\} dw}_{\text{integrates to 1}} \nonumber \\
& =
\frac{1}{\sqrt{\gamma}} \exp{\left\{- \frac{\gamma}{\mu} \right\} } \exp{\left\{
\frac{a(y-b)}{c} \right\}} \hspace{1in} (\text{for $\gamma,\mu>0$}),
\label{eq:pdfIntw}
\end{align}
%---------------------------
where the third line makes the substitutions $\gamma=\frac{a^{2} + 2c}{c
}=\frac{1}{2p(1-p)}$, $\mu^{2}=\frac{\gamma c}{(y-b)^{2}}=\frac{\sigma^{2}}
{p^{2}(1-p)^{2}(y-b)^{2}}$, and $\mu=\frac{\sigma}{p(1-p)|y-b|}$. In the
fourth line, integration with respect to $w$ yields 1 because it is the
\emph{pdf} of a reciprocal inverse-Gaussian distribution i.e. $w\sim
RIG(\gamma,\mu)$. Substituting the values of $(\gamma,\mu,a,b,c)$ in
equation~\eqref{eq:pdfIntw} and canceling terms we get,
%-----------------------------
\begin{align}
P &=
\sqrt{2p(1-p)} \; \exp{\left\{ \frac{(1-2p)(y-b)}{2\sigma}- \frac{|y-b|}{2\sigma} \right\}}
 	\nonumber \\
&= \sqrt{2p(1-p)} \; \exp{\left\{ -\frac{1}{\sigma}\left[ p-I(y \le b)
	\right](y-b)\right\}} \nonumber\\
&=
\sqrt{2p(1-p)} \; \exp{\left\{ -\frac{1}{\sigma}\left[ p-I(y\le\mu+\sigma\alpha s)
	\right](y-\mu-\sigma\alpha s)\right\}}. \label{eq:pdfPvalue}
\end{align}
%-----------------------------
Substituting the value of $P$ from equation~\eqref{eq:pdfPvalue} into
equation~\eqref{eq:pdfDblInt}, canceling terms, writing the \emph{pdf} of S
and letting $\kappa = 2p(1-p)/\sigma$, the \emph{pdf} of $Y$ is,
%----------------------------
\begin{align}
f(y|\theta)& =
%-------------
\kappa \int_{\mathbb{R^+}} \exp\left\{ -\frac{1}{\sigma}
\left[p-I(y\le\mu+\sigma\alpha s) \right] (y-\mu-\sigma\alpha s)\right\}
\frac{1}{\sqrt{2\pi}} \, \exp\left\{-\frac{s^2}{2}\right\} ds
\nonumber \\
& = \kappa \int_0^\infty \frac{1}{\sqrt{2 \pi}}
\exp\left\{ -\frac{1}{2} \bigg[ s^2 + 2\left(\frac{y-\mu}{\sigma}
-\alpha s \right)\left[p-I\Big(\frac{y-\mu}{\sigma} \le \alpha s \Big) \right]
\bigg] \right\} ds.
\label{eq:pdf4cases}
\end{align}
%----------------------------
Evaluation of the \emph{pdf} $f(y|\theta)$ given by
equation~\eqref{eq:pdf4cases} leads to 4 cases depending on the sign of
$\alpha$ and $y^{\ast}=(y-\mu)/\sigma$ and we integrate them one at a time.
We also employ the earlier introduced notation $p_{\alpha_{-}}=
p-I(\alpha<0)$ and $p_{\alpha_{+}}= p-I(\alpha>0)$ in each cases.

\textbf{Case (i)}: When $(\alpha>0, y^{\ast}\le0)$, then $I(y^{\ast} \le
\alpha s) = 1$. The corresponding \emph{pdf} is,
%----------------
\begin{align}
f(y|\theta)& =
\kappa \int_0^\infty \frac{1}{\sqrt{2\pi}}
\exp\left\{ -\frac{1}{2} \left[ s^2 + 2\left(\frac{y-\mu}{\sigma}-\alpha s
\right)(p-1) \right] \right\} ds  \nonumber \\
& =
\kappa \int_0^\infty \frac{1}{\sqrt{2\pi}}
\exp\left\{ -\frac{1}{2} \left[ s^2 + 2\left(y^{\ast}-\alpha s
\right)p_{\alpha_{+}} \right] \right\} ds \nonumber \\
& =
\kappa \exp\left\{-y^{\ast}p_{\alpha_{+}} +
\frac{1}{2} \alpha^2 p_{\alpha_{+}}^2 \right\}
\int_0^\infty \frac{1}{\sqrt{2\pi}}
\exp\left\{-\frac{1}{2} \left(s-\alpha p_{\alpha_{+}}\right)^2 \right\} ds
\nonumber \\
& =
\kappa \exp\left\{-y^{\ast}p_{\alpha_{+}} +
\frac{1}{2} \alpha^2 p_{\alpha_{+}}^2 \right\}
\Phi\left(s - \alpha p_{\alpha_{+}}\right)\Big|_{0}^{\infty} \nonumber \\
& =
\kappa \, \Phi\left(\alpha p_{\alpha_{+}} \right) \exp{ \left\{-y^{*}
p_{\alpha_{+}}+ \frac{1}{2} \alpha^{2}
p_{\alpha_{+}}^{2} \right\}}. \label{eq:pdfCase1}
\end{align}
%---------------------

\textbf{Case (ii)}: For $(\alpha>0, y^{\ast}>0)$ we have two cases. Case (a):
$y^{\ast} > \alpha s$ implies $I(y^{\ast}\le\alpha s)=0$ and this occurs for
all $s \in [0, y^{\ast}/\alpha)$. Case(b): $y^{\ast} \le \alpha s$ implies
$I(y^{\ast}\le \alpha s)=1$ and this occurs for all $s \in
[y^{\ast}/\alpha,\infty)$. Hence the \emph{pdf} is,
%------------------------------
\begin{align}
f(y|\theta) & =
\kappa \int_0^{\frac{y^*}{\alpha}}
\frac{1}{\sqrt{2\pi}}\exp\left\{ -\frac{1}{2} \left[ s^2 +
2\left(y^{\ast}-\alpha s \right)p \right] \right\} ds \; +
\nonumber \\
& \hspace{17pt} \kappa \int_{\frac{y^*}{\alpha} }^{\infty} \frac{1}{\sqrt{2\pi}}
\exp\left\{ -\frac{1}{2} \left[ s^2 + 2\left(y^{\ast}-\alpha s
\right)(p-1) \right] \right\} ds \nonumber \\
%----		
& = \kappa  \exp\{-y^{*} p_{\alpha_{-}}\} \int_0^{\frac{y^*}{\alpha}}
\frac{1}{\sqrt{2\pi}} \exp\left\{ -\frac{1}{2} \left(s^2 -
2 \alpha p_{\alpha_{-}} s \right)\right\}ds \; + \nonumber \\
%-----	
& \hspace{17pt} \kappa \exp\{-y^{*}p_{\alpha_{+}}\}\int_{\frac{y^*}{\alpha} }^{\infty}
\frac{1}{\sqrt{2\pi}} \exp\left\{ -\frac{1}{2}
\left( s^2 - 2\alpha p_{\alpha_{+}} s\right)\right\} ds \nonumber \\
%-----
& = \kappa \exp\left\{-y^{*}p_{\alpha_{-}} + \frac{1}{2} \alpha^{2}
p_{\alpha_{-}}^{2} \right\} \int_0^{\frac{y^*}{\alpha}} \frac{1}{\sqrt{2\pi}}
\exp\left\{ -\frac{1}{2} \left(s-\alpha p_{\alpha_{-}} \right)^{2}\right\}
ds \; +\nonumber \\
%----
& \hspace{17pt} \kappa \exp\left\{-y^{*}p_{\alpha_{+}} + \frac{1}{2} \alpha^2
p_{\alpha_{+}}^{2} \right\} \int_{\frac{y^*}{\alpha}}^{\infty}
\frac{1}{\sqrt{2\pi}}\exp\left\{ -\frac{1}{2}	
\left( s^2 -\alpha p_{\alpha_{+}} \right)^2\right\} ds \nonumber \\
%----	
& = \kappa \exp\left\{-y^{*} p_{\alpha_{-}} + \frac{1}{2} \alpha^{2}
p_{\alpha_{-}}^{2} \right\} \Phi\left(s-\alpha p_{\alpha_{-}} \right)
\Big|_{0}^{\frac{y^{\ast}}{\alpha}} \; + \nonumber\\
%----
& \hspace{17pt}	\kappa \exp\left\{-y^{*}p_{\alpha_{+}} + \frac{1}{2}
\alpha^{2} p_{\alpha_{+}}^{2}\right\}
\Phi\left(s-\alpha p_{\alpha_{+}} \right)
\big|_{\frac{y^{\ast}}{\alpha}}^{\infty} \nonumber \\
%----	
& = \kappa \left[\Phi\left(\frac{y^*}{\alpha}- \alpha p_{\alpha_{-}} \right)
-\Phi(-\alpha p_{\alpha_{-}}) \right]
\exp\left\{-y^{*}p_{\alpha_{-}} + \frac{1}{2}
\alpha^{2} p_{\alpha_{-}}^{2}\right\} \; + \notag \\
%----	
&  \hspace{17pt} \kappa \left[ \Phi\left(\alpha p_{\alpha_{+}} -\frac{y^*}{\alpha}\right)
\right]	\exp\left\{-y^{*}p_{\alpha_{+}} + \frac{1}{2}
\alpha^{2} p_{\alpha_{+}}^{2} \right\}.\label{eq:pdfCase2}
\end{align}
%----------------------------

\textbf{Case (iii)}: When $(\alpha<0, y^{\ast}>0)$, then $I(y^{\ast}
\le\alpha s) = 0$ since $\alpha<0$. Hence we have,
%----------------------------
\begin{align}
f(y|\theta) & =
\kappa \int_{0}^{\infty} \frac{1}{\sqrt{2 \pi}} \exp\left\{ -\frac{1}{2}
\left[s^{2} + 2(y^{\ast}-\alpha s)p_{\alpha_{+}}  \right] \right\} ds \nonumber\\
& =
\kappa \exp\left\{-y^{\ast} p_{\alpha_{+}} + \frac{1}{2} \alpha^{2}
p_{\alpha_{+}}^{2}\right\} \int_{0}^{\infty} \frac{1}{\sqrt{2 \pi}}
\exp\left\{ -\frac{1}{2} \left(s-\alpha p_{\alpha_{+}} \right)^{2} \right\} ds
\nonumber \\
& =
\kappa \, \Phi(\alpha p_{\alpha_{+}}) \exp\left\{-y^{\ast} p_{\alpha_{+}}
+ \frac{1}{2} \alpha^{2} p_{\alpha_{+}}^{2}\right\}. \label{eq:pdfCase3}
\end{align}

%------------------------------------------------------------------------------
\textbf{Case (iv)}: For $(\alpha<0, y^{\ast}\le 0)$ we have two cases. Case
(a): $y^{\ast} \le \alpha s$ implies $I(y^{\ast} \le \alpha s)=1$ and this
occurs for all $s \in [0, y^{\ast}/\alpha]$. Case(b): $y^{\ast} > \alpha s$
implies $I(y^{\ast} \le \alpha s)=0$ and this occurs for all $s \in
(y^{\ast}/\alpha,\infty)$. So we have,
%----------------------------
\begin{align}
f(y|\theta) & =
\kappa \int_{0}^{\frac{y^{\ast}}{\alpha}} \frac{1}{\sqrt{2 \pi}}
\exp\left\{-\frac{1}{2} \left[ s^{2} + 2(y^{\ast} - \alpha s)(p-1) \right]
\right\} ds \; + \nonumber \\
& \hspace{17pt} \kappa \int_{\frac{y^{\ast}}{\alpha}}^{\infty} \frac{1}{\sqrt{2 \pi}}
\exp\left\{-\frac{1}{2} \left[ s^{2} + 2(y^{\ast} - \alpha s)p \right]
\right\} ds \nonumber \\
%----
& = \kappa \int_{0}^{\frac{y^{\ast}}{\alpha}} \frac{1}{\sqrt{2 \pi}}
\exp\left\{-\frac{1}{2} \left[ s^{2} + 2(y^{\ast} - \alpha s)p_{\alpha_{-}}
\right] \right\} ds \; + \nonumber \\
& \hspace{17pt} \kappa \int_{\frac{y^{\ast}}{\alpha}}^{\infty} \frac{1}{\sqrt{2 \pi}}
\exp\left\{-\frac{1}{2} \left[ s^{2} + 2(y^{\ast} - \alpha s) p_{\alpha_{+}}
\right] \right\} ds \nonumber \\
%----	
& = \kappa \left[\Phi\left(\frac{y^*}{\alpha}- \alpha p_{\alpha_{-}} \right)
-\Phi(-\alpha p_{\alpha_{-}}) \right]
\exp\left\{-y^{*}p_{\alpha_{-}} + \frac{1}{2}
\alpha^{2} p_{\alpha_{-}}^{2}\right\} \; + \nonumber \\
%----	
& \hspace{17pt} \kappa \left[ \Phi\left(\alpha p_{\alpha_{+}} -\frac{y^*}{\alpha}\right)
\right]	\exp\left\{-y^{*}p_{\alpha_{+}} + \frac{1}{2}
\alpha^{2} p_{\alpha_{+}}^{2} \right\},
\label{eq:pdfCase4}
\end{align}
%----------------------------
where the integration details are similar to Case~(ii) and have been
suppressed to avoid monotonicity and save space.

Combining all the four cases, i.e. equations \eqref{eq:pdfCase1} to
\eqref{eq:pdfCase4}, we have the \emph{pdf} of the GAL distribution given by
equation~\eqref{eq:GALpdf}. \hfill $\blacksquare$

%------------------------------------------------------------------------------
%------------------------------------------------------------------------------
\subsection{Cumulative Distribution Function}\label{app:GALcdf}

\textbf{Theorem 2}: Suppose $Y \sim GAL(\mu,\sigma,p,\alpha)$ and let
$y^{\ast}=(y-\mu)/\sigma$, then the \emph{cdf} $F$ is,
%---------------------------------
\begin{equation}
\begin{split}
F(y|\theta) = \bigg(1 - 2 \Phi\left(-\frac{y^{\ast}}{|\alpha|} \right)
- \frac{2p(1-p)}{p_{\alpha_{-}}} \exp\left\{- y^{\ast} p_{\alpha_{-}}
+ \frac{1}{2} \alpha^{2} p_{\alpha_{-}}^{2}  \right\}
\Big[\Phi\left(\frac{y^{\ast}}{\alpha} - \alpha p_{\alpha_{-}} \right) \\
%----
- \, \Phi(-\alpha p_{\alpha_{-}}) \Big] \bigg) I\left( \frac{y^{\ast}}{\alpha}>0
\right) + I(\alpha <0) - \frac{2p(1-p)}{p_{\alpha_{+}}}
\exp\left\{-y^{\ast} p_{\alpha_{+}} + \frac{1}{2} \alpha^{2} p_{\alpha_{+}}^{2}
\right\} \\
%-----
\times \, \Phi\left[\alpha p_{\alpha_{+}} - \frac{y^{\ast}}{\alpha}
I\left( \frac{y^{\ast}}{\alpha}>0\right)   \right].
\end{split} \label{eq:cdf}
\end{equation}
%---------------------------------

\textbf{Proof}: We note that for any \emph{cdf} $F(y|\theta) =
\int\limits_{-\infty}^{y} f(v|\theta) dv = 1 - \int\limits_{y}^{\infty}
f(v|\theta) dv$. This property is used in deriving the \emph{cdf} when
$y>\mu$ to avoid breaking the region of integration as $(-\infty,\mu) \bigcup
(\mu,y)$. We let $v^{\ast}=(v-\mu)/\sigma$, combine cases and derive as
follows.

\textbf{Case (i)}: When ($\alpha>0,y\le \mu$) or ($\alpha<0, y>\mu$), the
\emph{cdf} is,
%----------------------------------
\begin{align}
F(y|\theta) & =
\begin{cases}
\int\limits_{-\infty}^{y} \kappa \Phi(\alpha p_{\alpha_{+}})
\exp\left\{ - v^{\ast} p_{\alpha_{+}} + \frac{1}{2} \alpha^{2}
p_{\alpha_{+}}^{2} \right\} dv, & \text{if $\alpha >0, y\le \mu$} \\
%----
1 - \int\limits_{y}^{\infty} \kappa \Phi(\alpha p_{\alpha_{+}})
\exp\left\{ - v^{\ast} p_{\alpha_{+}} + \frac{1}{2} \alpha^{2}
p_{\alpha_{+}}^{2} \right\} dv, & \text{if $\alpha <0, y > \mu$} \nonumber
\end{cases} \\
& =
\begin{cases}
\kappa \Phi(\alpha p_{\alpha_{+}})
\exp\left\{ \frac{1}{2} \alpha^{2} p_{\alpha_{+}}^{2} \right\}
\left[ \frac{\exp\left\{- v^{\ast} p_{\alpha_{+}} \right\}}
{- p_{\alpha_{+}}/\sigma} \right]_{-\infty}^{y},
& \text{if $\alpha >0, y \le \mu$}\\
%----
1 - \kappa \Phi(\alpha p_{\alpha_{+}}) \exp\left\{ \frac{1}{2} \alpha^{2}
p_{\alpha_{+}}^{2} \right\} \left[ \frac{\exp\left\{- v^{\ast} p_{\alpha_{+}} \right\}}
{- p_{\alpha_{+}}/\sigma} \right]_{y}^{\infty},
& \text{if $\alpha <0, y > \mu$}    \nonumber \\
\end{cases}\\
& =
\begin{cases}
\frac{2p(1-p)}{-p_{\alpha_{+}}} \,
\Phi(\alpha p_{\alpha_{+}}) \exp\left\{- y^{\ast} p_{\alpha_{+}}
+ \frac{1}{2} \alpha^{2} p_{\alpha_{+}}^{2} \right\},
& \text{if $\alpha >0, y \le \mu$}\\
%----
1 -  \frac{2p(1-p)}{p_{\alpha_{+}}} \,
\Phi(\alpha p_{\alpha_{+}}) \exp\left\{-y^{\ast} p_{\alpha_{+}}
+ \frac{1}{2} \alpha^{2} p_{\alpha_{+}}^{2} \right\},
& \text{if $\alpha <0, y > \mu$}
\end{cases}
\label{eq:cdfCase1}
\end{align}
%----------------------------------
where the third step substitutes the value of $\kappa$, $p_{\alpha_{+}} =
p-1$ for $\alpha>0$ and $p_{\alpha_{+}} = p$ for $\alpha<0$.

\textbf{Case (ii)}: When ($\alpha<0,y\le \mu$) or ($\alpha>0, y>\mu$), the
\emph{cdf} is,
%----------------------------------
\begin{align}
F(y|\theta) & =
\begin{cases}
\int\limits_{-\infty}^{y} \kappa \bigg( \left[ \Phi\left(\frac{v^{\ast}}
{\alpha} - \alpha p_{\alpha_{-}} \right) - \Phi(-\alpha p_{\alpha_{-}}) \right]
\exp\left\{-v^{\ast} p_{\alpha_{-}} + \frac{1}{2} \alpha^{2} p_{\alpha_{-}}^{2}
\right\} \\
+ \; \Phi\left(\alpha p_{\alpha_{+}} - \frac{v^{\ast}}{\alpha} \right)
\exp\left\{-v^{\ast} p_{\alpha_{+}} + \frac{1}{2} \alpha^{2}
p_{\alpha_{+}}^{2} \right\} \bigg) dv,
\hspace{0.5in} \text{if $\alpha<0, y\le\mu$}  \\
%-------------------------
1- \int\limits_{y}^{\infty} \kappa \bigg( \left[ \Phi\left(\frac{v^{\ast}}
{\alpha} - \alpha p_{\alpha_{-}} \right) - \Phi(-\alpha p_{\alpha_{-}}) \right]
\exp\left\{-v^{\ast} p_{\alpha_{-}} + \frac{1}{2} \alpha^{2} p_{\alpha_{-}}^{2}
\right\} \\
+ \; \Phi\left(\alpha p_{\alpha_{+}} - \frac{v^{\ast}}{\alpha} \right)
\exp\left\{-v^{\ast} p_{\alpha_{+}} + \frac{1}{2} \alpha^{2}
p_{\alpha_{+}}^{2} \right\} \bigg) dv,
\hspace{0.5in} \text{if $\alpha>0, y>\mu$}. \label{eq:cdfCase2Eval}
\end{cases}
\end{align}
%----------------------------------
Note that in both the subcases of equation~\eqref{eq:cdfCase2Eval}, the
integral remains the same and only the limits of integration changes. Hence,
we evaluate each terms individually over the limits ($a,b$) and will
substitute values of $(a,b)$ as per our requirement.

To evaluate the first integral component of equation~\eqref{eq:cdfCase2Eval}
denoted $C_{1}$, we substitute $z = v^{\ast}/\alpha  - \alpha p_{\alpha_{-}}$
and perform integration-by-parts as follows.
%--------------------------
\begin{align}
C_{1} & = \kappa \int\limits_{a}^{b} \Phi\left(\frac{v^{\ast}}
{\alpha} - \alpha p_{\alpha_{-}}  \right) \exp\left\{ -v^{\ast}
p_{\alpha_{-}} + \frac{1}{2} \alpha^{2} p_{\alpha_{-}}^{2}   \right\} dv
\nonumber\\
& = \kappa \exp\left\{-\frac{1}{2} \alpha^{2} p_{\alpha_{-}}^{2} \right\}
\int_{\frac{a^{\ast}}{\alpha} - \alpha p_{\alpha_{-}}}^{\frac{b^{\ast}}{\alpha}
- \alpha p_{\alpha_{-}}} \left[ \alpha \sigma \Phi(z) \exp\{-\alpha
p_{\alpha_{-}} z\} \right] dz
\nonumber \\
& = \kappa \sigma \exp\left\{-\frac{1}{2} \alpha^{2} p_{\alpha_{-}}^{2}
\right\} \Bigg( \bigg[\frac{\alpha \Phi(z) \exp\left\{-\alpha p_{\alpha_{-}} z
\right\} }{-\alpha p_{\alpha_{-}}}  \bigg]_{\frac{a^{\ast}}{\alpha}
- \alpha p_{\alpha_{-}}}^{\frac{b^{\ast}}{\alpha} - \alpha p_{\alpha_{-}}}
\nonumber \\
& \hspace{14pt} - \, \frac{\alpha}{-\alpha p_{\alpha_{-}}}
\int_{\frac{a^{\ast}}{\alpha} - \alpha p_{\alpha_{-}}}^{\frac{b^{\ast}}
{\alpha} - \alpha p_{\alpha_{-}}} \frac{1}{\sqrt{2\pi}} \exp\left\{-\frac{1}{2}
(z^{2} + 2\alpha p_{\alpha_{-}} z + \alpha^{2} p_{\alpha_{-}}^{2}
- \alpha^{2} p_{\alpha_{-}}^{2})  \right\} dz
\Bigg) \nonumber \\
& = \kappa \sigma \exp\left\{-\frac{1}{2}
\alpha^{2} p_{\alpha_{-}}^{2} \right\} \Bigg( \frac{\exp\left\{ \alpha^{2}
p_{\alpha_{-}}^{2} \right\}}{- p_{\alpha_{-}}} \bigg[\Phi\left(\frac{b^{\ast}}{\alpha} -
\alpha p_{\alpha_{-}} \right) \exp\left\{-b^{\ast} p_{\alpha_{-}} \right\}
\nonumber \\
& \hspace{14pt}
- \Phi\left(\frac{a^{\ast}}{\alpha}
- \alpha p_{\alpha_{-}} \right) \exp\left\{-a^{\ast} p_{\alpha_{-}} \right\}
\bigg] + \frac{1}{p_{\alpha_{-}}} \exp\left\{\frac{1}{2}
\alpha^{2} p_{\alpha_{-}}^{2} \right\}
\Phi\left(z + \alpha p_{\alpha_{-}} \right)_{\frac{a^{\ast}}{\alpha} -
\alpha p_{\alpha_{-}}}^{\frac{b^{\ast}}{\alpha} - \alpha p_{\alpha_{-}}}
\Bigg)
\nonumber \\
& = -\frac{ \kappa \sigma}{p_{\alpha_{-}}} \exp\left\{\frac{1}{2}
\alpha^{2} p_{\alpha_{-}}^{2} \right\} \bigg[\Phi\left(\frac{b^{\ast}}
{\alpha} - \alpha p_{\alpha_{-}} \right) \exp\left\{-b^{\ast}
p_{\alpha_{-}} \right\}
\nonumber \\
& \hspace{14pt}
- \Phi\left(\frac{a^{\ast}}{\alpha}
- \alpha p_{\alpha_{-}} \right) \exp\left\{-a^{\ast} p_{\alpha_{-}} \right\}
\bigg] + \frac{\kappa \sigma}{p_{\alpha_{-}}}
\left[\Phi\left(\frac{b^{\ast}}{\alpha}\right) -
\Phi\left(\frac{a^{\ast}}{\alpha} \right) \right].
\label{eq:cdfC1}
\end{align}
%--------------------------
We next evaluate the second component denoted $C_{2}$ directly as follows,
%--------------------------
\begin{align}
C_{2} & = - \kappa \Phi(-\alpha p_{\alpha_{-}}) \int_{a}^{b}  \exp\left\{
-v^{\ast} p_{\alpha_{-}} + \frac{1}{2} \alpha^{2} p_{\alpha_{-}}^{2}\right\} dv
\nonumber \\
& = - \kappa \Phi(-\alpha p_{\alpha_{-}}) \exp\left\{ \frac{1}{2} \alpha^{2}
p_{\alpha_{-}}^{2}\right\}  \left[\frac{\exp\left\{-v^{\ast} p_{\alpha_{-}}
\right\} }{-p_{\alpha_{-}}/\sigma}   \right]_{a}^{b}
\nonumber \\
&=  \frac{\kappa \sigma}{p_{\alpha_{-}}} \Phi(-\alpha p_{\alpha_{-}})
\exp\left\{ \frac{1}{2} \alpha^{2} p_{\alpha_{-}}^{2} \right\}
\left[ \exp\left\{-b^{\ast} p_{\alpha_{-}} \right\} -
\exp\left\{-a^{\ast} p_{\alpha_{-}} \right\} \right].
\label{eq:cdfC2}
\end{align}
%--------------------------
Finally, to evaluate the third component denoted $C_{3}$ we use the
substitution $z=\alpha p_{\alpha_{+}} - v^{\ast}/\alpha$ and
integrate-by-parts as done in $C_{1}$. We suppress the details for brevity
and present the final expression:
%--------------------------
\begin{align}
C_{3} & = \kappa \int_{a}^{b} \Phi\left(\alpha p_{\alpha_{+}} -
\frac{v^{\ast}}{\alpha}  \exp\left\{-v^{\ast} p_{\alpha_{+}}
+ \frac{1}{2} \alpha^{2} p_{\alpha_{+}}^{2}  \right\}  \right) dv
\nonumber \\
& = - \frac{\kappa \sigma}{p_{\alpha_{+}}} \exp\left\{\frac{1}{2}
\alpha^{2} p_{\alpha_{+}}^{2} \right\} \left[ \Phi\left(\alpha p_{\alpha_{+}}
- \frac{b^{\ast}}{\alpha} \right) \exp\{-b^{\ast} p_{\alpha_{+}} \}
- \Phi\left(\alpha p_{\alpha_{+}} - \frac{a^{\ast}}{\alpha} \right)
\exp\{-a^{\ast} p_{\alpha_{+}} \} \right]
\nonumber \\
& \hspace{14pt} + \; \frac{\kappa \sigma}{p_{\alpha_{+}}} \left[ \Phi\left(-\frac{b^{\ast}}
{\alpha} \right) - \Phi\left(-\frac{a^{\ast}}{\alpha} \right) \right].
\label{eq:cdfC3}
\end{align}
%--------------------------
Note that
$\Phi\left(\frac{b^{\ast}}{\alpha}\right)-\Phi\left(\frac{a^{\ast}}{\alpha}
\right) = -\Phi\left(-\frac{b^{\ast}}{\alpha}\right)+
\Phi\left(-\frac{a^{\ast}}{\alpha}\right)$ and hence the relevant term from
equation~\eqref{eq:cdfC1} and equation~\eqref{eq:cdfC3} can be collected
together when adding the expressions.

When $\alpha<0$ and $y\le \mu$, the limits of integration $a=-\infty$ and
$b=y$ implies $a^{\ast}=-\infty$ and $b^{\ast}=y^{\ast}$, respectively.
Substituting the values of $a^{\ast}, b^{\ast}$ and $\kappa$ in $C_{1},
C_{2}$ and $C_{3}$ and summing the expression yields,
%--------------------------
\begin{align}
F(y|\theta) & =
%-----
2\Phi\left(-\frac{y^{\ast}}{\alpha} \right) - \frac{2p(1-p)}{p_{\alpha_{-}}}
\exp\left\{-y^{\ast}p_{\alpha_{-}} + \frac{1}{2} \alpha^{2} p_{\alpha_{-}}^{2}
\right\}  \left[\Phi\left(\frac{y^{\ast}}{\alpha} - \alpha p_{\alpha_{-}} \right)
 - \Phi(-\alpha p_{\alpha_{-}})  \right] \nonumber \\
%-----
& \hspace{14pt} -\frac{2p(1-p)}{p_{\alpha_{+}}} \exp\left\{-y^{\ast} p_{\alpha_{+}} +
\frac{1}{2} \alpha^{2} p_{\alpha_{+}}^{2}   \right\}
\Phi\left(\alpha p_{\alpha_{+}} -\frac{y^{\ast}}{\alpha}  \right).
 \label{eq:cdfCase2i}
\end{align}
%--------------------------
Similarly, when $\alpha>0$ and $y > \mu$, the limits of integration $a=y$ and
$b=\infty$ implies $a^{\ast}=y^{\ast}$ and $b^{\ast}=\infty$, respectively.
Substituting the values of $a^{\ast}, b^{\ast}$ and $\kappa$ in $C_{1},
C_{2}$ and $C_{3}$ and evaluating the expression $1-C_{1}-C_{2}-C_{3}$,
yields,
%--------------------------
\begin{align}
F(y|\theta) & =
%-----
1-2\Phi\left(-\frac{y^{\ast}}{\alpha} \right) - \frac{2p(1-p)}{p_{\alpha_{-}}}
\exp\left\{-y^{\ast}p_{\alpha_{-}} + \frac{1}{2} \alpha^{2} p_{\alpha_{-}}^{2}
\right\}  \bigg[\Phi\left(\frac{y^{\ast}}{\alpha} - \alpha p_{\alpha_{-}} \right)
 \nonumber \\
%-----
& \hspace{14pt}  - \, \Phi(-\alpha p_{\alpha_{-}})  \bigg]
-\frac{2p(1-p)}{p_{\alpha_{+}}} \exp\left\{-y^{\ast} p_{\alpha_{+}} +
\frac{1}{2} \alpha^{2} p_{\alpha_{+}}^{2}   \right\}
\Phi\left(\alpha p_{\alpha_{+}} -\frac{y^{\ast}}{\alpha}  \right).
\label{eq:cdfCase2ii}
\end{align}
%--------------------------
Combining the equations \eqref{eq:cdfCase1}, \eqref{eq:cdfCase2i} and
\eqref{eq:cdfCase2ii}, we have the \emph{cdf} of the GAL distribution given
by equation~\eqref{eq:GALcdf}. \hfill $\blacksquare$

%------------------------------------------------------------------------------
%------------------------------------------------------------------------------
\subsection{Moment Generating Function}\label{app:GALmgf}

\textbf{Theorem 3}: Suppose $Y \sim GAL(\mu,\sigma,p,\alpha)$, then the
\emph{mgf} denoted $M_{Y}(t)$ is as follows,
%---------------------------------
\begin{equation}
M_{Y}(t) = 2p(1-p) \bigg[\frac{ (p_{\alpha_{+}} - p_{\alpha_{-}}) }
{(p_{\alpha_{-}} -\sigma t)(p_{\alpha_{+}} -\sigma t)} \bigg]
\exp\Big\{ \mu t + \frac{1}{2} \alpha^2 \sigma^2 t^2 \Big\} \;
\Phi\big(|\alpha| \sigma t \big). \label{eq:mgf}
\end{equation}
%---------------------------------

\textbf{Proof}: Using the definition of the \emph{mgf} we have,
%---------------------------------
\begin{equation}
M_{Y}(t) = \int_{-\infty}^{\infty} \exp(ty)
f(y|\mu,\sigma,p,\alpha) \, dy.
\label{eq:mgfdef}
\end{equation}
%---------------------------------
Substituting the GAL \emph{pdf} \eqref{eq:GALpdf} into
equation~\eqref{eq:mgfdef}, leads to two cases depending on $\alpha > 0$ or
$\alpha < 0$.  We again use the notation $\kappa = 2p(1-p)/\sigma$ and break
the region of integration depending on $y^{\ast} > 0$ (i.e., $y>\mu$) or
$y^{\ast} \le 0$ (i.e., $y\le\mu$).

\textbf{Case (i)}: When $\alpha > 0$ and $y^{\ast} > 0$ (i.e., $y>\mu$), we
have the following three components,
%---------------------------------
\begin{align}
M_{1} & = \kappa \int_{\mu}^{\infty} \exp\{ty\} \,
\Phi\Big(\frac{y^{\ast}}{\alpha} -
\alpha p_{\alpha_{-}} \Big) \exp\left\{-y^{\ast} p_{\alpha_{-}} +
\frac{1}{2} (\alpha p_{\alpha_{-}} )^{2}   \right\} dy,  \nonumber \\
M_{2} & = - \kappa  \int_{\mu}^{\infty} \exp\{ty\} \,
\Phi(-\alpha p_{\alpha_{-}})
\exp\left\{-y^{\ast} p_{\alpha_{-}} + \frac{1}{2}(\alpha p_{\alpha_{-}} )^{2}
\right\} dy, \label{eq:M1-M3Int}   \\
M_{3} & = \kappa  \int_{\mu}^{\infty} \exp\{ty\} \,
\Phi\bigg( \alpha p_{\alpha_{+}}
- \frac{y^{\ast}}{\alpha} \bigg) \exp\left\{-y^{\ast} p_{\alpha_{+}} +
\frac{1}{2}(\alpha p_{\alpha_{+}} )^{2}   \right\} dy, \nonumber
\end{align}
%---------------------------------
and when $\alpha > 0$ and $y^{\ast} \le 0$ (i.e., $y \le \mu$) we have,
\begin{equation}
M_{4} =  \kappa  \int_{-\infty}^{\mu}  \exp\{ty\} \,
\Phi\big( \alpha p_{\alpha_{+}} \big)
\exp\left\{ -y^{\ast} p_{\alpha_{+}} +
\frac{1}{2}(\alpha p_{\alpha_{+}} )^{2}   \right\} dy. \label{eq:M4Int}
\end{equation}
We first consider $M_{1}$, substitute $z = \frac{y^{\ast}}{\alpha} - \alpha
p_{\alpha_{-}}$, change the limits of integration and integrate-by-parts as
follows:
%----------------------------------
\begin{align}
 M_{1}  & =  \int_{-\alpha p_{\alpha_{-}}}^{\infty} \kappa \alpha \sigma \,
\exp\big\{ \mu t + \alpha \sigma (z + \alpha p_{\alpha_{-}}) t \big\}
\Phi(z) \exp\left\{ - \alpha p_{\alpha_{-}} (z + \alpha p_{\alpha_{-}})
+ \frac{1}{2} \alpha^{2} p_{\alpha_{-}}^{2} \right\} dz \nonumber \\
& = \kappa \alpha \sigma \exp\bigg\{ \mu t + \alpha^{2} p_{\alpha_{-}}
\sigma t - \frac{1}{2}\alpha^{2} p_{\alpha_{-}}^{2} \bigg\}
\int_{-\alpha p_{\alpha_{-}}}^{\infty} \Phi(z)
\exp\left\{-\alpha(p_{\alpha_{-}} - \sigma t)z \right\}  dz \nonumber \\
& = \kappa \alpha \sigma \exp\bigg\{ \mu t + \alpha^{2} p_{\alpha_{-}}
\sigma t - \frac{1}{2}\alpha^{2} p_{\alpha_{-}}^{2} \bigg\}
\Bigg[\Phi(z) \, \frac{\exp\{-\alpha(p_{\alpha_{-}} - \sigma t)z\}}
{-\alpha(p_{\alpha_{-}} - \sigma t)}\bigg|_{-\alpha p_{\alpha_{-}}}^{\infty}
\nonumber \\
& \hspace{14pt} - \int_{-\alpha p_{\alpha_{-}}}^{\infty} \phi(z) \;
\frac{\exp\{-\alpha(p_{\alpha_{-}} - \sigma t)z\}}
{-\alpha(p_{\alpha_{-}} - \sigma t)} \; dz \Bigg]
\hspace{1in} \text{(limit exists only if $t<p/\sigma$)}\nonumber \\
& = \kappa \sigma \frac{\exp\big\{ \mu t + \alpha^{2} p_{\alpha_{-}}
\sigma t - \frac{1}{2}\alpha^{2} p_{\alpha_{-}}^{2} \big\}}{(p_{\alpha_{-}}
- \sigma t)} \Bigg[\exp\{\alpha^{2}p_{\alpha_{-}} (p_{\alpha_{-}} - \sigma t)\}
\Phi(-\alpha p_{\alpha_{-}})  \; +  \nonumber \\
& \hspace{14pt} \exp\bigg\{\frac{\alpha^{2}}{2} (p_{\alpha_{-}} - \sigma t)^{2}\bigg\}
\int_{-\alpha p_{\alpha_{-}}}^{\infty} \frac{1}{\sqrt{2 \pi}}
\exp\bigg\{-\frac{1}{2} \big(z+\alpha(p_{\alpha_{-}} - \sigma t)^{2} \big)
\bigg\} \; dz \Bigg] \nonumber\\
& = \kappa \sigma \frac{\exp\big\{ \mu t + \frac{1}{2}\alpha^{2}
p_{\alpha_{-}}^{2} \big\}}{(p_{\alpha_{-}} - \sigma t)}
\Phi(-\alpha p_{\alpha_{-}}) +
\kappa \sigma \frac{\exp\big\{ \mu t + \frac{1}{2}
\alpha^{2} \sigma^{2} t^{2} \big\}}{(p_{\alpha_{-}} - \sigma t)} \;
\Phi\big(z+\alpha(p_{\alpha_{-}} - \sigma t)\big)
\Big|_{-\alpha p_{\alpha_{-}}}^{\infty} \nonumber \\
& = \kappa \sigma \frac{\exp\big\{ \mu t + \frac{1}{2}\alpha^{2}
p_{\alpha_{-}}^{2} \big\}}{(p_{\alpha_{-}} - \sigma t)} \;
\Phi(-\alpha p_{\alpha_{-}})
+ \kappa \sigma
\frac{\exp\big\{ \mu t + \frac{1}{2}\alpha^{2} \sigma^{2} t^{2}
\big\}}{(p_{\alpha_{-}} - \sigma t)} \; \Phi(\alpha \sigma t). \label{eq:M1}
\end{align}

Secondly, we integrate the expression for $M_{2}$ as follows,
%---------------------------------
\begin{align}
M_{2} & = - \kappa \int_{\mu}^{\infty} \exp\{ty\} \,
\Phi(-\alpha p_{\alpha_{-}})
\exp\left\{-y^{\ast} p_{\alpha_{-}} + \frac{1}{2} \alpha^{2} p_{\alpha_{-}}^{2}
\right\} dy \nonumber \\
& = - \kappa \; \Phi(-\alpha p_{\alpha_{-}}) \exp\bigg\{
\frac{p_{\alpha_{-}} \mu}{\sigma} + \frac{1}{2} \alpha^{2} p_{\alpha_{-}}^{2}\bigg\}
\int_{\mu}^{\infty} \exp\left\{ - \frac{(p_{\alpha_{-}} - \sigma t)}{\sigma} \; y
\right\}  dy \nonumber \\
& = - \kappa \; \Phi(-\alpha p_{\alpha_{-}}) \exp\bigg\{
\frac{p_{\alpha_{-}} \mu}{\sigma} + \frac{1}{2} \alpha^{2} p_{\alpha_{-}}^{2}\bigg\}
\times \frac{(-1)(- \sigma)}{(p_{\alpha_{-}} - \sigma t)} \exp\bigg\{
\frac{-(p_{\alpha_{-}} - \sigma t)}{\sigma} \mu \bigg\} \nonumber \\
& = - \kappa \sigma \frac{ \exp\left\{ \mu t + \frac{1}{2} \alpha^{2}
p_{\alpha_{-}}^{2} \right\}}{(p_{\alpha_{-}}
- \sigma t)}\;  \Phi(-\alpha p_{\alpha_{-}})
\hspace{0.5in} \text{(limit exists only if $t<p/\sigma$)}.
\label{eq:M2}
\end{align}
%---------------------------------

The integral $M_{3}$ is evaluated by substituting $z = \alpha p_{\alpha_{+}}
- \frac{y^{\ast}}{\alpha}$. Thereafter, the integration is analogous to
$M_{1}$ with the existence condition $t>(p-1)/\alpha$ and yields,
%---------------------------------
\begin{equation}
M_{3} = \kappa \sigma \frac{\exp\big\{ \mu t + \frac{1}{2}\alpha^{2}
p_{\alpha_{+}}^{2} \big\}}{(p_{\alpha_{+}} - \sigma t)} \;
\Phi(\alpha p_{\alpha_{+}})
- \kappa \sigma \frac{\exp\big\{ \mu t + \frac{1}{2}\alpha^{2}
\sigma^{2} t^{2} \big\}}{(p_{\alpha_{+}} - \sigma t)} \; \Phi(\alpha \sigma t).
\label{eq:M3}
\end{equation}
%---------------------------------
The integral $M_{4}$ can be directly evaluated and closely follows the steps
in the integration of $M_{2}$ with the resulting expression,
%---------------------------------
\begin{equation}
M_{4} = - \kappa \sigma \frac{ \exp\left\{ \mu t + \frac{1}{2}
\alpha^{2} p_{\alpha_{+}}^{2} \right\}}{(p_{\alpha_{+}}
- \sigma t)} \; \Phi(\alpha p_{\alpha_{+}}). \label{eq:M4}
\end{equation}
%---------------------------------
where the limits of integration exists only for $t>(p-1)/\sigma$.

To obtain the \emph{mgf} of the GAL distribution, we substitute the values of
$\kappa$ and sum the four equations \eqref{eq:M1}, \eqref{eq:M2},
\eqref{eq:M3} and \eqref{eq:M4}. This yields,
%---------------------------------
\begin{equation}
M_{Y}(t) = 2p(1-p) \bigg[ \frac{(p_{\alpha_{+}} -
p_{\alpha_{-}})}{(p_{\alpha_{-}} - \sigma t)(p_{\alpha_{+}}-\sigma t)} \bigg]
\exp\bigg\{ \mu t +\frac{1}{2} \alpha^{2}\sigma^{2} t^{2}
\bigg\} \; \Phi(\alpha \sigma t) .
\label{eq:mgfalphaPos}
\end{equation}
%---------------------------------

\textbf{Case (ii)}: When $\alpha < 0$ and $y^{\ast}\le0$ (i.e., $y\le\mu$),
we have the following three components of the \emph{mgf},
%---------------------------------
\begin{align}
M_{5} & = \kappa \int_{-\infty}^{\mu} \exp\{ty\} \,
\Phi\Big(\frac{y^{\ast}}{\alpha} -
\alpha p_{\alpha_{-}} \Big) \exp\left\{-y^{\ast} p_{\alpha_{-}} +
\frac{1}{2} (\alpha p_{\alpha_{-}} )^{2}   \right\} dy,  \nonumber \\
M_{6} & = - \kappa \int_{-\infty}^{\mu}  \exp\{ty\} \,
\Phi(-\alpha p_{\alpha_{-}})
\exp\left\{-y^{\ast} p_{\alpha_{-}} + \frac{1}{2}(\alpha p_{\alpha_{-}} )^{2}
\right\} dy, \label{eq:M5-M7Int}   \\
M_{7} & = \kappa \int_{-\infty}^{\mu} \exp\{ty\} \,
\Phi\bigg( \alpha p_{\alpha_{+}}
- \frac{y^{\ast}}{\alpha} \bigg) \exp\left\{-y^{\ast} p_{\alpha_{+}} +
\frac{1}{2}(\alpha p_{\alpha_{+}} )^{2}   \right\} dy, \nonumber
\end{align}
%---------------------------------
and when $\alpha < 0$ and $y^{\ast}>0$ (i.e., $y>\mu$) we have,
\begin{equation}
M_{8} =  \kappa \int_{\mu}^{\infty}  \exp\{ty\} \,
\Phi\big( \alpha p_{\alpha_{+}} \big)
\exp\left\{ -y^{\ast} p_{\alpha_{+}} +
\frac{1}{2}(\alpha p_{\alpha_{+}} )^{2}   \right\}  dy. \label{eq:M8Int}
\end{equation}
%---------------------------------
The integration of the terms $M_{5}$ to $M_{8}$ are similar to the case when
$\alpha>0$ and results in the following \emph{mgf},
%---------------------------------
\begin{equation}
M_{Y}(t) = 2p(1-p) \bigg[ \frac{(p_{\alpha_{+}} -
p_{\alpha_{-}})}{(p_{\alpha_{-}} - \sigma t)(p_{\alpha_{+}}-\sigma t)} \bigg]
\exp\bigg\{ \mu t +\frac{1}{2} \alpha^{2}\sigma^{2} t^{2}
\bigg\} \; \Phi(-\alpha \sigma t) .
\label{eq:mgfalphaNeg}
\end{equation}
%---------------------------------

Combining the \emph{mgf} for the two cases, i.e. equations
\eqref{eq:mgfalphaPos} and \eqref{eq:mgfalphaNeg}, we have the \emph{mgf} of
the GAL distribution. \hfill $\blacksquare$

%------------------------------------------------------------------------------
%------------------------------------------------------------------------------
The mean, variance and skewness of the GAL distribution can be obtained from
the GAL \emph{mgf} \eqref{eq:mgf}. We state this in terms of a theorem below.

\textbf{Theorem 4}: Suppose $Y \sim GAL(\mu,\sigma,p,\alpha)$, then
%--------------------------
\begin{equation}
\begin{split}
E(Y) & = \mu + \frac{2|\alpha| \sigma}{\sqrt{2 \pi}} + \sigma
         \frac{(1-2p)}{p(1-p)}    \\
V(Y) & = \alpha^{2} \sigma^{2} \bigg(1- \frac{2}{\pi} \bigg) + \sigma^{2}
         \bigg[\frac{1-2p + 2p^{2}}{p^{2}(1-p)^{2}}  \bigg]            \\
S(Y) & = \frac{\alpha^{3}  \sqrt{\frac{2}{\pi}} \left(\frac{4}{\pi} - 1 \right)
+ 2 \left[ \frac{(1-p)^{3}-p^{3}}{p^{3}(1-p)^{3}} \right] }
{ \left\{ \alpha^{2} \left(1- \frac{2}{\pi} \right) +
\left[ \frac{1-2p+2p^{2}}{p^{2}(1-p)^{2}} \right] \right\}^{3/2}}
\end{split}
\label{Theorem4}
\end{equation}
%--------------------------
where $E(Y)$, $V(Y)$ and $S(Y)$ denote the mean, variance and skewness,
respectively.

\textbf{Proof}: Taking logarithm of the \emph{mgf} and keeping terms
involving $t$ we have,
%--------------------------
\begin{equation}
\ln M_{Y} (t) \propto -\ln(p_{\alpha_{-}} - \sigma t) - \ln(p_{\alpha_{+}} - \sigma t)
+ \mu t + \frac{1}{2}\alpha^{2} \sigma^{2} t^{2} + \ln
\Phi(|\alpha| \sigma t).
\label{eq:mgfprop}
\end{equation}
%--------------------------
Taking the first, second and third derivative of equation~\eqref{eq:mgfprop}
and evaluating at $t=0$ we get,
%--------------------------
\begin{alignat*}{3}
m_{1}(Y) & = E(Y) = & \frac{\partial \ln M_{Y}(t)}{\partial t} \Big|_{t=0}
& =  \mu + \sqrt{\frac{2}{\pi}} \; \alpha \sigma + \sigma \frac{(1-2p)}{p(1-p)} \\
m_{2}(Y) & = V(Y) = & \frac{\partial^{2} \ln M_{Y}(t)}{\partial t^{2}}\Big|_{t=0}
& =  \alpha^{2} \sigma^{2} \left(1-\frac{2}{\pi} \right) + \sigma^{2}
\left[\frac{1-2p+2p^{2}}{p^{2}(1-p)^{2}} \right] \\
m_{3}(Y) & =  & \frac{\partial^{3} \ln M_{Y}(t)}{\partial t^{3}}\Big|_{t=0} & =
\alpha^{3} \sigma^{3}\sqrt{\frac{2}{\pi}} \left(\frac{4}{\pi} - 1 \right)
+ 2 \sigma^{3} \left[ \frac{(1-p)^{3} - p^{3} }{p^{3}(1-p)^{3}}   \right]
\end{alignat*}
%--------------------------
where $m_{1}(Y), m_{2}(Y)$ and $m_{3}(Y)$ are the first, second and third
order central moments, respectively. Hence, skewness can be obtained as
$\frac{m_{3}}{(m_{2}^{3/2})}$. \hfill $\blacksquare$

%------------------------------------------------------------------------------
%------------------------------------------------------------------------------
\section{\textbf{Conditional Densities in the FBQROR model}}\label{app:FBQROR}

In this appendix, we derive the conditional posteriors of the FBQROR model
parameters. Specifically, the conditional posteriors of $\beta,\nu,h$, and
$z$ have tractable distributions and is sampled using a Gibbs approach. The
parameters $(\sigma,\gamma)$ are jointly sampled using random-walk MH
algorithm (to reduce autocorrelation in MCMC draws) and $\delta$ is sampled
using a random-walk MH algorithm. The derivations below follow the ordering
as presented in Algorithm~\ref{alg:algorithm1}.

%-------------------------------------------
\textbf{(1)} Starting with $\beta$, the conditional posterior
$\pi(\beta|z,\nu,h,\sigma,\gamma)$ is proportional to $\pi(\beta) \times
f(z|\beta,\nu,h,\sigma,\gamma)$ and its kernel can be written as,
%--------------------------
\begin{align*}
\pi(\beta|z,\nu,h,\sigma,\gamma) & \propto
\exp\bigg\{-\frac{1}{2} \bigg[ \sum_{i=1}^{n} \frac{(z_{i}-x'_{i}\beta -
A\nu_{i}- C |\gamma| h_{i})^2}{\sigma B \nu_{i}}
+ (\beta - \beta_{0})'B_{0}^{-1}(\beta - \beta_{0})\bigg]  \bigg\} \\
%-----
& \propto \exp \bigg\{ -\frac{1}{2} \bigg[ \beta'\bigg(\sum_{i=1}^{n}
\frac{x_{i}x'_{i}}{\sigma B \nu_{i}} + B_{0}^{-1} \bigg) \beta
-  \beta'  \bigg( \frac{x_{i}(z_{i}-A\nu_{i}- C |\gamma| h_{i})}
{\sigma B \nu_{i}}  + B_{0}^{-1}\beta_{0} \bigg) \\
& \hspace{44pt}  - \bigg( \frac{x'_{i}(z_{i}-A\nu_{i}-C |\gamma| h_{i})}
{\sigma B \nu_{i}}  + \beta'_{0}B_{0}^{-1} \bigg)
\bigg]   \bigg\} \\
%-----
& \propto \exp \bigg\{-\frac{1}{2} \bigg[\beta' \tilde{B}^{-1}\beta -
\beta'\tilde{B}^{-1} \tilde{\beta} - \tilde{\beta}'\tilde{B}^{-1}\beta
+ \tilde{\beta}'\tilde{B}^{-1}\tilde{\beta}
- \tilde{\beta}'\tilde{B}^{-1}\tilde{\beta} \bigg] \bigg\}\\
%-----
& \propto \exp \Big\{-\frac{1}{2} (\beta-\tilde{\beta})'\tilde{B}^{-1}
(\beta-\tilde{\beta})  \Big\},
\end{align*}
%--------------------------
where the posterior variance $\tilde{B}$ and the posterior mean
$\tilde{\beta}$ are defined as follows:
%--------------------------
\begin{equation*}
\tilde{B}^{-1} = \bigg(\sum_{i=1}^{n} \frac{x_{i} x'_{i}}{\sigma
         B \nu_{i}} + B_{0}^{-1} \bigg)
         \hspace{0.25in} \mathrm{and} \hspace{0.25in}
         \tilde{\beta} =
         \tilde{B}\bigg( \sum_{i=1}^{n} \frac{x_{i}(z_{i}
         - A \nu_{i} - C|\gamma| h_{i})}{\sigma B \nu_{i}}
         + B_{0}^{-1} \beta_{0} \bigg).
\end{equation*}
%--------------------------
Hence, the conditional posterior is a normal distribution and
$\beta|z,\nu,h,\sigma,\gamma \sim N(\tilde{\beta}, \tilde{B})$.

%-------------------------------------------
\textbf{(2)} The parameters $(\sigma,\gamma)$ are jointly sampled marginally
of $(z,\nu,h)$ to reduce autocorrelation in the MCMC draws. Collecting terms
involving $(\sigma,\gamma)$ from the complete data posterior
density~\eqref{eq:compDataPosterior} does not yield a tractable distribution,
hence $(\sigma, \gamma)$ are sampled using a joint random-walk MH algorithm.
The proposed values are generated from a truncated bivariate normal
distribution $TBN_{(0,\infty) \times (L,U)}\big( (\sigma_{c},\gamma_{c}),
\iota_{1}^{2} \hat{D}_{1}\big)$, where $(\sigma_{c},\gamma_{c})$ represent
the current values, $\iota_{1}$ denotes the tuning factor and $\hat{D}_{1}$
is the negative inverse of the Hessian obtained by maximizing the
log-likelihood~\eqref{eq:fulllikelihood} with respect to $(\sigma, \gamma)$.
This maximization process to obtain $\hat{D}_{1}$ is computationally
expensive and so may be done only once at the beginning of the algorithm with
$\beta$ values fixed at the BQROR or OLS estimates. The proposed draws are
accepted with MH probability,
%--------------------------
\begin{equation*}
\alpha_{MH}(\sigma_{c},\gamma_{c};\sigma',\gamma') = \min
\bigg\{0,\ln\bigg[\frac{f(y| \beta,\sigma',\gamma',\delta) \,
\pi(\beta,\sigma',\gamma',\delta)}{f(y|
\beta,\sigma_{c},\gamma_{c},\delta) \, \pi(\beta,\sigma_{c},\gamma_{c},\delta)} \;
\frac{\pi(\sigma_{c},\gamma_{c}| (\sigma',\gamma'), \iota_{1}^{2} \hat{D}_{1})}
{\pi(\sigma',\gamma'| (\sigma_{c},\gamma_{c}), \iota_{1}^{2} \hat{D}_{1})}
\bigg] \bigg\},
\end{equation*}
%--------------------------
otherwise, $(\sigma_{c},\gamma_{c})$ is repeated in the next MCMC iteration.
Here, $f(\cdot)$ represents the full likelihood \eqref{eq:fulllikelihood},
$\pi(\beta,\sigma,\delta,\gamma)$ denotes the prior
distributions~\eqref{eq:priors}, $\pi(\sigma_{c},\gamma_{c}|
(\sigma',\gamma'), \iota_{1}^{2} \hat{D}_{1})$ stands for the truncated
bivariate normal probability with mean $(\sigma',\gamma')$ and covariance
$\iota_{1}^{2} \hat{D}_{1}$. The expression $\pi(\sigma',\gamma'|
(\sigma_{c},\gamma_{c}), \iota_{1}^{2} \hat{D}_{1})$ has an analogous
interpretation. Note that the tuning parameter $\iota_{1}$ can be adjusted
for appropriate step-size and acceptance rate and the parameters $(A,B,C)$
depend on $p$ which in turn is a function of $p_{0}$ and $\gamma$.

%-------------------------------------------
\textbf{(3)} The  conditional posterior of $\nu$ is obtained from the
complete posterior density~\eqref{eq:compDataPosterior} by collecting terms
involving $\nu$. This is done element-wise as follows:
%--------------------------
\begin{align}
\pi(\nu_{i}|z,\beta,h,\sigma,\gamma)
& \propto \nu_{i}^{-\frac{1}{2}} \exp \bigg\{-\frac{1}{2} \bigg[
\frac{(z_{i}-x'_{i}\beta - A\nu_{i}- C|\gamma|h_{i})^{2}}
{\sigma B \nu_{i}}  \bigg] -\frac{\nu_{i}}{\sigma}  \bigg\} \notag \\
%-----
& \propto \nu_{i}^{-\frac{1}{2}} \exp \bigg\{-\frac{1}{2} \bigg[
\frac{(z_{i}-x'_{i}\beta - C|\gamma|h_{i})^{2}}
{\sigma B } \, \nu_{i}^{-1} + \bigg(\frac{A^{2}}{\sigma B}
+ \frac{2}{\sigma} \bigg) \nu_{i} \bigg] \bigg\} \notag \\
%-----
& \propto \nu_{i}^{-\frac{1}{2}} \exp\bigg\{-\frac{1}{2}
\bigg[a_{i} \nu_{i}^{-1} + b \nu_{i}  \bigg]  \bigg\}, \notag
\end{align}
%--------------------------
which is recognized as the kernel of a  generalized inverse-Gaussian (GIG)
distribution where,
%--------------------------
\begin{equation*}
a_{i} = \frac{(z_{i}-x'_{i}\beta - C|\gamma|h_{i})^{2}}{\sigma B}
         \hspace{0.25in} \mathrm{and} \hspace{0.25in}
b = \bigg(\frac{A^{2}}{\sigma B} + \frac{2}{\sigma} \bigg).
\end{equation*}
%--------------------------
Hence, $\nu_{i}|z,\beta,h,\sigma,\gamma \sim GIG(0.5,a_{i},b)$ for
$i=1,2,\ldots,n$.

%-------------------------------------------
\textbf{(4)} The  conditional posterior of $h$ is obtained element-wise from
the complete posterior density \eqref{eq:compDataPosterior} conditional on
$\gamma$ and remaining parameters as follows:
%--------------------------
\begin{align}
\pi(h_{i}|z,\beta,\nu,\sigma,\gamma)
& \propto \exp \bigg\{ -\frac{1}{2} \bigg[
\frac{(z_{i}-x'_{i}\beta - A\nu_{i}- C|\gamma|h_{i})^{2}}
{\sigma B \nu_{i}}  + \frac{h_{i}^{2}}{\sigma^2} \bigg]\bigg\} \notag \\
%-----
& \propto \exp \bigg\{-\frac{1}{2} \bigg[ \bigg(\frac{1}{\sigma^2} +
\frac{C^{2} \gamma^{2}}{\sigma B \nu_{i}}  \bigg) h_{i}^{2} -
\frac{2  C|\gamma|(z_{i}-x'_{i}\beta -
A\nu_{i})}{\sigma B \nu_{i}} h_{i}\bigg]  \bigg\} \notag \\
%-----
& \propto \exp \bigg\{-\frac{1}{2} \bigg[(\sigma_{h_{i}}^{2})^{-1}
h_{i}^{2}  - 2 (\sigma_{h_{i}}^{2})^{-1} \mu_{h_{i}} h_{i}\bigg]
\bigg\} \notag \\
%-----
& \propto \exp \bigg\{ -\frac{1}{2} (\sigma_{h_{i}}^{2})^{-1}
(h_{i} - \mu_{h_{i}})^{2} \bigg\}, \notag
\end{align}
%--------------------------
where the third line introduces the following notations:
%--------------------------
\begin{equation*}
(\sigma_{h_{i}}^{2})^{-1} = \bigg(\frac{1}{\sigma^2} + \frac{C^{2}
\gamma^{2}}{\sigma  B \nu_{i}}  \bigg)
         \hspace{0.25in} \mathrm{and} \hspace{0.25in}
\mu_{h_{i}} = \sigma_{h_{i}}^{2} \bigg( \frac{ C|\gamma|(z_{i}- x'_{i}\beta
- A\nu_{i})}{\sigma B \nu_{i}}  \bigg),
\end{equation*}
%--------------------------
and the fourth line adds and subtracts $(\sigma_{h_{i}}^{2})^{-1}
\mu_{h_{i}}^{2}$ to complete the square. The last expression is recognized as
the kernel of a half-normal distribution and hence,
$h_{i}|z,\beta,\nu,\sigma,\gamma \sim N^{+}(\mu_{h_{i}},\sigma_{h_{i}}^{2})$
for $i=1,2,\ldots,n$.

%-------------------------------------------
\textbf{(5)} The transformed cut-point $\delta$ is sampled from the full
likelihood \eqref{eq:fulllikelihood}, marginally of ($z,\nu,h$). The proposed
values are generated from a random-walk chain, $\delta'=\delta_{c} + u$,
where $u\sim N(0_{J-3}, \iota_{2}^{2}\hat{D}_{2})$, $\iota_{2}$ is a tuning
parameter and $\hat{D}_{2}$ denotes negative inverse Hessian, obtained by
maximizing the log-likelihood with respect to $\delta$. Given the current
value $\delta_{c}$, the proposed value $\delta'$ is accepted with MH
probability,
%--------------------------
\begin{equation*}
\alpha_{MH}(\delta_{c},\delta') = \min \bigg\{0,\ln\bigg[\frac{f(y|
\beta,\sigma,\gamma,\delta') \, \pi(\beta,\sigma,\gamma,\delta')}{f(y|
\beta,\sigma,\gamma,\delta_{c}) \, \pi(\beta,\sigma,\gamma,\delta_{c})}  \bigg]
\bigg\},
\end{equation*}
%--------------------------
otherwise, the current value $\delta_{c}$ is repeated. The variance of $u$
may be tuned as required for an appropriate step-size and acceptance rate.

%-------------------------------------------
\textbf{(6)} The full conditional density of the latent variable $z$ is a
truncated normal distribution where the cut-point vector $\xi$ is obtained
based on one-to-one mapping with $\delta$. Hence, $z$ is sampled as
$z_{i}|\beta,\nu,h,\sigma,\gamma,\delta,y \sim TN_{(\xi_{j-1},\,
\xi_{j})}(x'_{i}\beta+A \nu_{i}+  C|\gamma|h_{i}, \sigma B \nu_{i})$ for
$i=1,\cdots,n$ and $j=1,\cdots,J$.

%---------------------------------- Bibliography ------------------------------
\clearpage \pagebreak %\nocite{*}
\pdfbookmark[1]{References}{unnumbered} % To make References as a bookmark in pdf
\section*{References}

\bibliography{BibFBQROR}
\bibliographystyle{jasa}

%\bibliography{C:/1-ARSHAD-OPTIPLEX/Google-Drive/18-References/BibEconometrics}
%\bibliographystyle{C:/1-ARSHAD-OPTIPLEX/Google-Drive/18-References/jasa}

\end{document}